\documentclass[10pt,twoside,english,a4paper]{article}
\usepackage{graphicx} 
\usepackage{caption}
\usepackage{subcaption}
\usepackage[dvipsnames]{xcolor}
\usepackage{amsmath}
\usepackage{amsfonts}
\usepackage{amssymb}
\usepackage{mathtools}
\usepackage[edges]{forest} 
\usepackage{float} 
\usepackage{gensymb} 

\usepackage{booktabs}

\usepackage{comment}
\usepackage{placeins}
\usepackage[T1]{fontenc}
\usepackage{textcomp}
\usepackage{gensymb}

\usepackage{booktabs}
\usepackage{siunitx}
\usepackage{threeparttable}

\newcommand{\freqheader}{%
  \toprule
  {Reference} &
  {$\rho_{\partial \Omega}$} &
  \multicolumn{10}{c}{$\omega$ [M\,rad/s]} &
  {Average [\%]} \\
  \cmidrule(lr){3-12}
  {} & {[kg/$\rm m^2$]} &
  {1.26} & {2.51} & {3.77} & {5.03} & {6.28} &
  {7.54} & {8.80} & {10.05} & {11.31} & {12.57} & {} \\
  \midrule
}

\usepackage[english]{babel}
\usepackage{csquotes}

\DeclareQuoteStyle{englishguillemets}
{\guillemotleft}{\guillemotright}
{\guilsinglleft}{\guilsinglright}

\setquotestyle{englishguillemets}

\usepackage[
	commentmarkup=todo,
    todonotes={
        textwidth=1.75cm,
        textsize=scriptsize
    },
	commandnameprefix=always,
	defaultcolor=RedOrange
]{changes}
\definechangesauthor[name={Leo}, color=RedOrange]{LP}

\usepackage{soul} 

\sethighlightmarkup{%
  \IfIsColored
    {{\sethlcolor{authorcolor!30}\hl{#1}}}
    {#1}%
}

\usepackage{fancyhdr}
\usepackage{authblk}

\usepackage[
  backend=biber,
  style=chem-acs,
  sorting=ynt,
  sortcites=true
]{biblatex}
\newcommand{\fig}[1]{Figure \ref{#1}} 
\newcommand{\equ}[1]{equation (\ref{#1})} 
\newcommand{\rmm}{relaxed micromorphic model } 
\newcommand{\rrmm}{reduced relaxed micromorphic model } 

\newcommand{\dyniprod}[1]{\left\langle #1 \right\rangle}

\newcommand{\Cmacro}{{\mathbb{C}_{\rm macro}}}
\newcommand{\Cmicro}{{\mathbb{C}_{\rm micro}}}
\newcommand{\evalat}[2]{{\left. #1 \right|_{#2}}}
\newcommand{\surfdens}{{\rho_{\partial \Omega}}}

\title{Breaking Scale Separation in Metamaterials' homogenization: Interface-Inertia-Enhanced Relaxed Micromorphic Model}

\author[,1]{Angela Madeo\thanks{Corresponding author. Email: \texttt{angela.madeo@tu-dortmund.de}}}

\author[1]{Leonardo Perez}

\author[1]{Mohammad Sarhil}

\affil[1]{Institute for Structural Mechanics and Dynamics, Faculty of Architecture and Civil Engineering, 
TU Dortmund, Germany}

\usepackage{subfiles} 

\begin{document}
	\maketitle
	
	\abstract{
		Homogenized continuum models are widely used to describe wave propagation and band-gap behavior in mechanical metamaterials without explicitly resolving their microstructure.
		Their validity, however, typically relies on the classical separation of scales assumption, according to which the wavelength of the propagating disturbance is much larger than the characteristic size of the unit cell.
		In finite-size metamaterial samples and at higher frequencies, this assumption progressively breaks down, and the dynamic response becomes strongly influenced by the way the microstructure is truncated at the external boundaries.
		
		In this work we introduce a fundamentally new concept in the homogenized description of mechanical metamaterials: the inertial contribution of macroscopic interfaces.
		We show that different truncations of the same lattice generate boundaries with distinct mass distributions, which lead to measurable differences in the dynamic response of finite-sized specimens.
		To capture this complex mechanism in a homogenized framework, we extend the relaxed micromorphic model by introducing a kinetic surface energy defined on the boundary of the considered body.
		This generates an additional inertial term in the boundary conditions that can be seen as the homogenized counterpart of the interface inertia produced by the truncation of the microstructure.
		As a result, the homogenized model can now distinguish between finite-sized specimens that share identical bulk properties but differ only in the configuration of their interfaces.
		
		The proposed formulation preserves the variational structure of the relaxed micromorphic model while enabling the continuum to reproduce boundary-dependent responses observed in fully resolved simulations, particularly in frequency regimes where scale separation no longer holds.
		
		This brand-new result opens the door to predictive homogenized models, enabling the upscale towards larger scale structures assembled from finite metamaterial building blocks, a step that would be computationally infeasible using fully resolved microstructured models.
	}
	
\section*{Keywords}

Interface inertia, Homogenization, Mechanical metamaterials,   Relaxed micromorphic model,  Boundary effects,  Finite-size metamaterials
	
\FloatBarrier

\section{Introduction}
	
	Mechanical metamaterials derive their unusual dynamic properties from the interaction between microstructural geometry and elastic wave propagation   \cite{Chaplain_2025}. Their microstructures can be tailored to guide, attenuate, and focus waves for targeted frequencies \cite{DorKanDreKoc:2026:gpm,KanDorDreKoc:2026:mam,PanBraCorSanDal:2026:tma,PanBraCorSaDal:2026:bsa,HerChMem:2026:tlo,DalVinMin:2026:hbs}. 
	In particular, band-gaps \, \textemdash frequency ranges in which waves cannot propagate\textemdash \, can emerge due to local resonance mechanisms \cite{Gao.2022,Sugino.2016,Sugino.2017,Sugino.2018,GeiBigMovBac:2024:bgp} or Bragg scattering effects \cite{Nobrega.2016,Krushynska.2017,Liu.2012,Sigalas.2005,Wen.2020,Zhang.2023,AleGou:2026:stc}. 
	 The homogenization of the dynamic response of the mechanical metamaterial is therefore an active field of research (see e.g. \cite{FraBigPic:2026:hoa,KucGeeKoz:2026:ijs,LiuVanGeeKou:2025:cma,Neff.2020,Demetriou.2025,Demore.2022}). 
	These phenomena have motivated the development of advanced continuum models  \cite{Min:1964:msi,MinEsh:1968:ofsg} capable of reproducing dispersive responses that cannot be captured by classical Cauchy elasticity, whose dispersion relations are strictly linear (i.e. non-dispersive).
	
	Among the available enriched continuum theories, the relaxed micromorphic model has proven particularly successful in describing dispersive wave propagation and band-gap behavior \cite{Sarhil.2026,Rizzi.2026,Demetriou.2025}.
	When calibrated appropriately, it reproduces the dispersion curves obtained from Bloch-Floquet analyses of periodic microstructured media and provides a computationally efficient homogenized description of the material response.
	Most validations of such models, however, rely on infinite or \enquote{effectively infinite} periodic domains, where the classical separation of scales assumption is satisfied: the characteristic size of the unit cell is assumed to be much smaller than the macroscopic dimensions of the body.
	
	In many practical applications of metamaterials, this assumption is violated.
	Finite samples are often composed of only a few unit cells, and their response becomes strongly influenced by the geometry of the external boundaries, particularly at higher frequencies where the wavelength approaches the characteristic length of the unit cell.
	In such situations, detailed numerical simulations that fully resolve the microstructure frequently reveal responses that differ substantially from those predicted by homogenized models calibrated on infinite periodic systems.
	These discrepancies become increasingly pronounced as the frequency increases and the classical scale-separation hypothesis progressively loses validity.
	
	A key aspect that remains largely unexplored in homogenized descriptions is the role played by the inertial contribution of the interface itself.
	Interface effects can be activated in finite-sized metamaterials when the wavelength of the propagating wave becomes sufficiently small to interact with the underlying microstructure \cite{Hermann.2026,PerezRamirez.2024}.
	In this regime, the wave no longer averages out the underlying geometry, and the details of how the lattice is truncated begin to significantly influence the global response.
	As a result, the interaction between the traveling wave and the microstructural length scale leads to boundary-driven phenomena that cannot be captured by classical homogenized models relying exclusively on bulk energies.
	
	In a lattice, truncating the periodic structure along different planes generates boundaries with markedly different mass distributions and geometric configurations.
	While these variations are naturally captured in fully resolved microstructured simulations, they are absent in standard continuum formulations, whose action functionals contain only bulk kinetic and strain energies.
	As a consequence, the influence of boundary inertia is not represented in the resulting equilibrium equations and boundary conditions.
	
	In this work we propose a homogenized relaxed micromorphic framework that explicitly incorporates the inertial contribution of solid interfaces.
	The central idea is to introduce a kinetic surface energy which can be put in relation with the boundary properties of the microstructured body, formulated within the variational structure of the macroscopic relaxed micromorphic model.
	This additional surface kinetic energy produces a boundary term in the equations of motion representing the effective inertia of the interface generated by truncating the underlying architected structure of the considered  macroscopic specimen. 
	The resulting formulation preserves the variational consistency of the model while enriching its boundary conditions with a physically motivated inertial contribution, allowing the relaxed micromorphic model to distinguish between specimens that share the same bulk material properties but differ solely in the way the lattice is truncated at their boundaries.
	
	Introducing interface inertia directly within the variational structure of a homogenized continuum model represents a fundamentally new concept in the homogenization of mechanical metamaterials.
	To the best of our knowledge, existing homogenized frameworks account exclusively for bulk energetic contributions, while boundary effects are treated only through classical traction conditions or external constraints.
	The present work is the first to embed an intrinsic inertial surface contribution within the Lagrangian of a homogenized medium, thereby extending the energetic structure of the continuum model itself to capture truncation-dependent, interface dynamic effects \footnote{We explicitly remark that many works exist in the literature that consider the effect of interface forces applied at the microscopic level (at interfaces inside the unit cell) on the macroscopic bulk behavior \cite{Murdoch.1976,Moeckel.1975,dellIsola.1987,Fried.2007,Gurtin.1975,ShaGan:2004:sde}. The approach that we use in this paper is quite different in the sense that here we consider interface forces to be applied at the macroscopic boundaries of the metamaterial’s specimen, so that they indeed represent a sort of homogenized characteristic of the interface itself.}.
	
	The framework developed here therefore offers a route to extend the applicability of homogenized  models to finite-sized metamaterial structures operating at higher frequencies, where classical scale separation is no longer strictly valid.
	By embedding interface inertia directly into the variational formulation, the proposed approach captures boundary-dependent dynamic responses while retaining the computational advantages of homogenized descriptions.
	
\FloatBarrier
	
\section{Description of the Problem}

	In the Bloch–Floquet analysis of a band-gap material, the dispersion curves are linear at sufficiently low frequencies and then progressively flatten as they approach the lower edge of the band gap.
	The relaxed micromorphic model has been shown to reproduce this dispersive behavior, both below and inside the band-gap frequency range \cite{Sarhil.2026,ErelDemore.2025,Demetriou.2025,PerezRamirez.2024,Demetriou.2024,PerezRamirez.2023}, especially when sufficiently large metamaterial specimens are considered \cite{Rizzi.2026,Demore.2022,Voss.2023}.
	In this regime, the classical separation of scales hypothesis is satisfied: the characteristic size of the unit cell is much smaller than the characteristic size of the macroscopic body and than the wavelength of the propagating disturbance.
	For finite-size metamaterial specimens, however, this assumption may no longer be valid.
	When the wavelength becomes comparable to the size of the specimen, and eventually to the size of the unit cell, the response of the structure can become sensitive to the details of the underlying microstructure.
	In this case, differences may arise between fully microstructured simulations and homogenized micromorphic predictions, especially close to the metamaterial’s boundaries.
	More importantly, specimens that share the same bulk periodic architecture may exhibit different dynamic responses if their external boundaries are generated by different cuts of the same lattice.

    In this work, we investigate this loss of scale separation by focusing on the role of the interface.
    We argue that one of the main sources of discrepancy in homogenized finite-size simulations of mechanical metamaterials is the inertial contribution associated with the boundary regions.
    This contribution is not accounted for in classical homogenized models, which are based only on bulk kinetic and strain energies.
    As a consequence, the boundary conditions obtained from the standard variational formulation cannot distinguish between specimens that have identical bulk properties but different interface configurations.
    
    \begin{figure}[!ht]
		\begin{subfigure}{\textwidth}
			\centering
			\includegraphics[width=\textwidth]{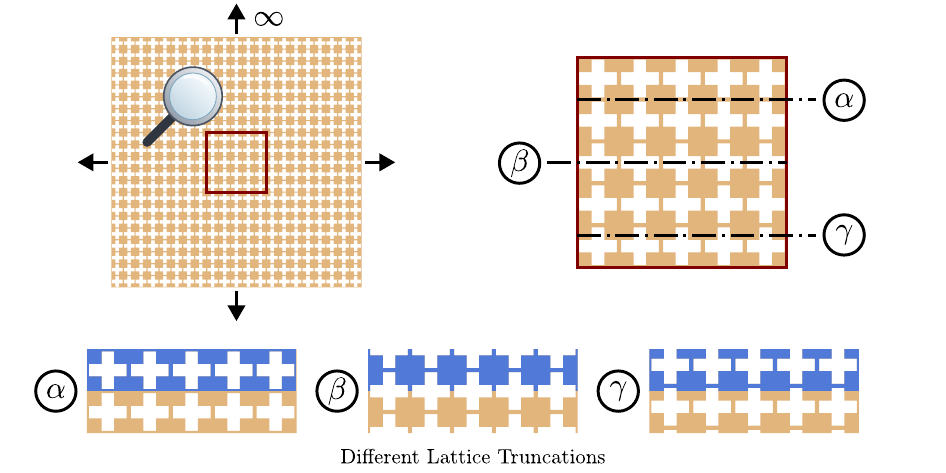}
		\end{subfigure}
		\caption{
			An infinite metamaterial is truncated using different planes; each boundary shows a different surface mass distribution.
			The figure also shows the three boundaries considered in this paper: $\alpha$-cut, $\beta$-cut, $\gamma-$cut.
		}
		\label{fig:interface_mass_distribution}
	\end{figure}

\begin{figure}[!ht]
		\begin{subfigure}{\textwidth}
			\centering
			\includegraphics[width=\textwidth]{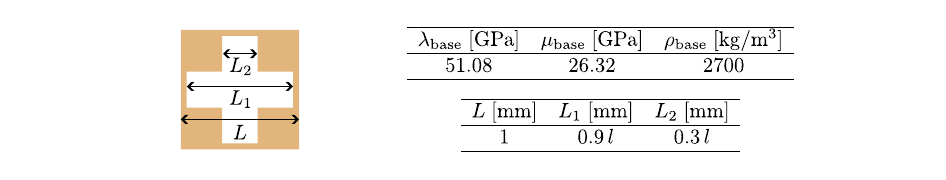}
		\end{subfigure}
		\caption{
			Dimension of the unit cell and geometric and material properties of the base material.
		}
		\label{fig:unit_cell-characteristics}
	\end{figure}

\begin{figure}[!ht]
		\begin{subfigure}{\textwidth}
			\centering
			\includegraphics[width= \textwidth]{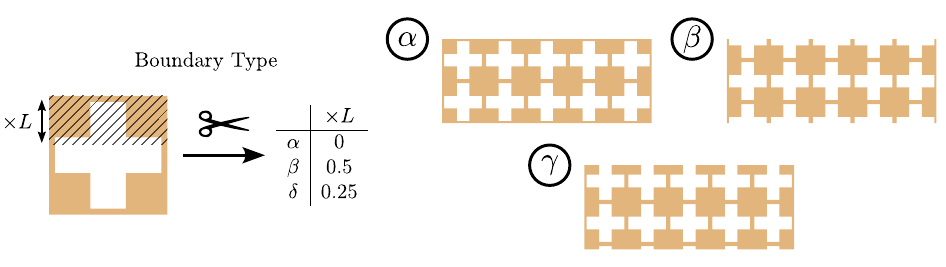}
		\end{subfigure}
		\caption{
			Three types of boundary, each associated with a different truncation plane.
		}
		\label{fig:unit_cell-boundaries}
	\end{figure}

    This situation naturally arises when an infinite periodic lattice is truncated along different planes (see \fig{fig:interface_mass_distribution}).
    Let us consider the metamaterial constructed from the unit cell shown in \fig{fig:interface_mass_distribution}, where the white regions represent voids distributed in a solid aluminum matrix (yellow).
    Depending on the position and orientation of the truncation plane, the resulting boundary can contain different amounts of solid material and therefore different mass distributions.
    Although these differences are automatically included in fully microstructured simulations, they are absent in a classical homogenized description.
    In the present paper we consider the "cross" unit cell studied in \cite{PerezRamirez.2023,Neff.2020,Rizzi.2022}, whose characteristics are given in \fig{fig:unit_cell-characteristics}.
    The metamaterial obtained by the periodic repetition of the "cross" unit cell can be truncated with different planes giving rise to different metamaterial interfaces.
    We consider three different metamaterial interfaces which we call $\alpha$, $\beta$, and $\gamma$, which are defined in \fig{fig:unit_cell-boundaries}.
    
    The benchmark problem considered here consists of a finite-sized metamaterial block embedded in a homogeneous isotropic medium and excited by an incident plane wave (see Figure \ref{fig:geometry_problem}).
    The interaction between the incoming wave and the finite microstructured block generates a scattered field, whose characteristics depend on both the bulk metamaterial properties and the geometry of the interface.
    Due to the discrete nature of the lattice, different cuts of the same periodic structure lead to different boundary configurations, so that different scattering patterns can be expected as soon as the separation of scales hypothesis breaks down.
    In the following subsection, we therefore compare fully microstructured simulations of two representative interface cuts, denoted as the $\alpha$-cut and the $\beta$-cut, in order to show how their dynamic responses differ as the frequency increases.
    Figure \ref{fig:unit_cell-characteristics} shows the geometric and material properties of the considered unit cell.
	
	\subsection{Influence of Boundary Truncation on the Scattering Response of finite-sized metamaterials}
	\label{sec:influence_boundary_truncation}
	
	In this subsection we present fully resolved microstructured simulations for a benchmark scattering problem in order to highlight on a concrete example how different interface cuts influence the dynamic response of finite metamaterial specimens.
	
    To this aim, we compare the results obtained from the two fully resolved simulations corresponding to the $\alpha$-cut and $\beta$-cut of the metamaterial.
    Both specimens share the same bulk microstructure and differ only in the way the lattice is truncated at the external boundary.
    
\begin{figure}[!ht]
		\begin{subfigure}{\textwidth}
			\centering
			\includegraphics[width=0.6 \textwidth]{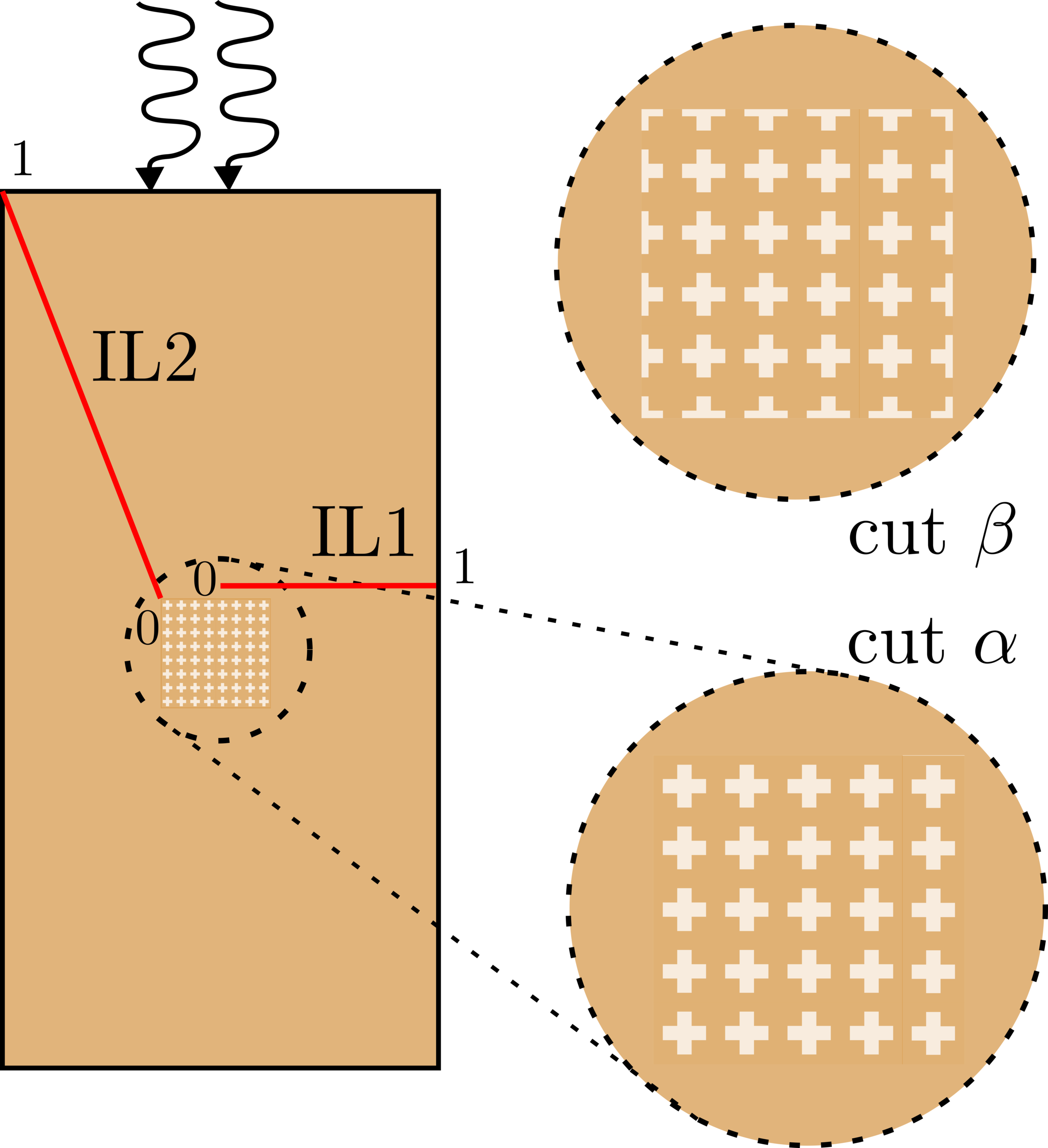}
		\end{subfigure}
		\caption{
			Definition of the benchmark problem studied in this paper. Inspection lines 1 and 2 are to give a quantitative  estimate of the normalized displacement of the $\alpha$- and $\beta-{\rm simulations}$ .
		}
		\label{fig:geometry_problem}
	\end{figure}

\begin{figure}[H]
	    \begin{subfigure}{\textwidth}
			\centering
			\includegraphics[width=\textwidth]{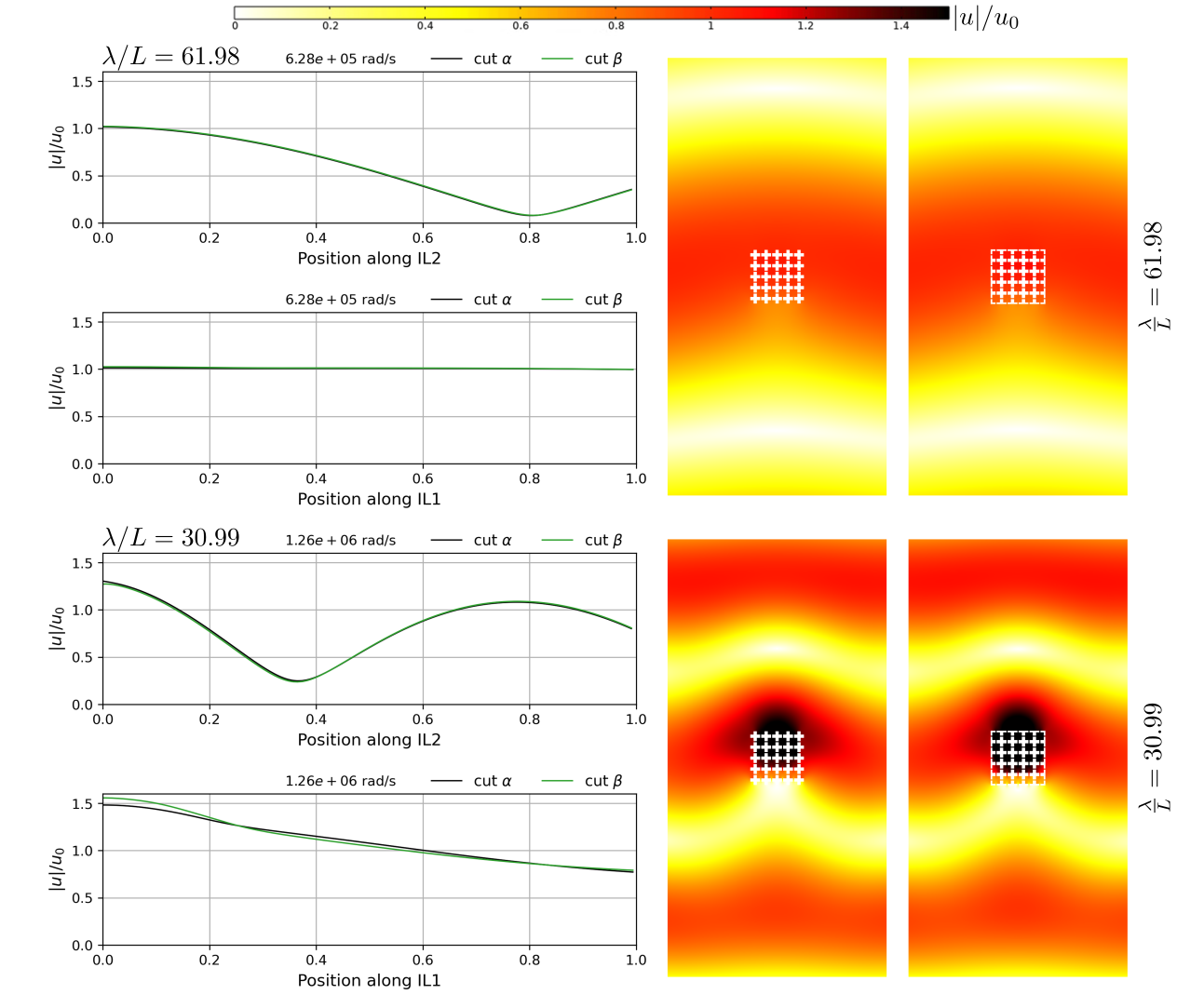}
		\end{subfigure}
		\caption{
		    Differentiated response of the scattering pattern, driven by the boundary contribution of $\alpha$ and $\beta$ interfaces, as a pressure wave interacts with a $5L\times5L$ metamaterial block.
		    The normalized displacement along inspection lines 1 and 2, as well as the normalized displacement field associated with a metamaterial block with $\alpha$- and $\beta$-type boundaries, are shown for the frequencies:
			\textit{(upper)} $0.63\ \mathrm{M\,rad/s}$ and
			\textit{(lower)} $1.26\ \mathrm{M\,rad/s}$.
			The corresponding wavelengths are much bigger than the size of the unit cell: no significant differences between $\alpha$- and $\beta$-cuts can be appreciated.
			The hypothesis of scale separation can be expected to hold true.
		}
		\label{fig:mstd_differences-5_low}
	\end{figure}

\begin{figure}[H]
		\begin{subfigure}{\textwidth}
			\centering
			\includegraphics[width=\textwidth]{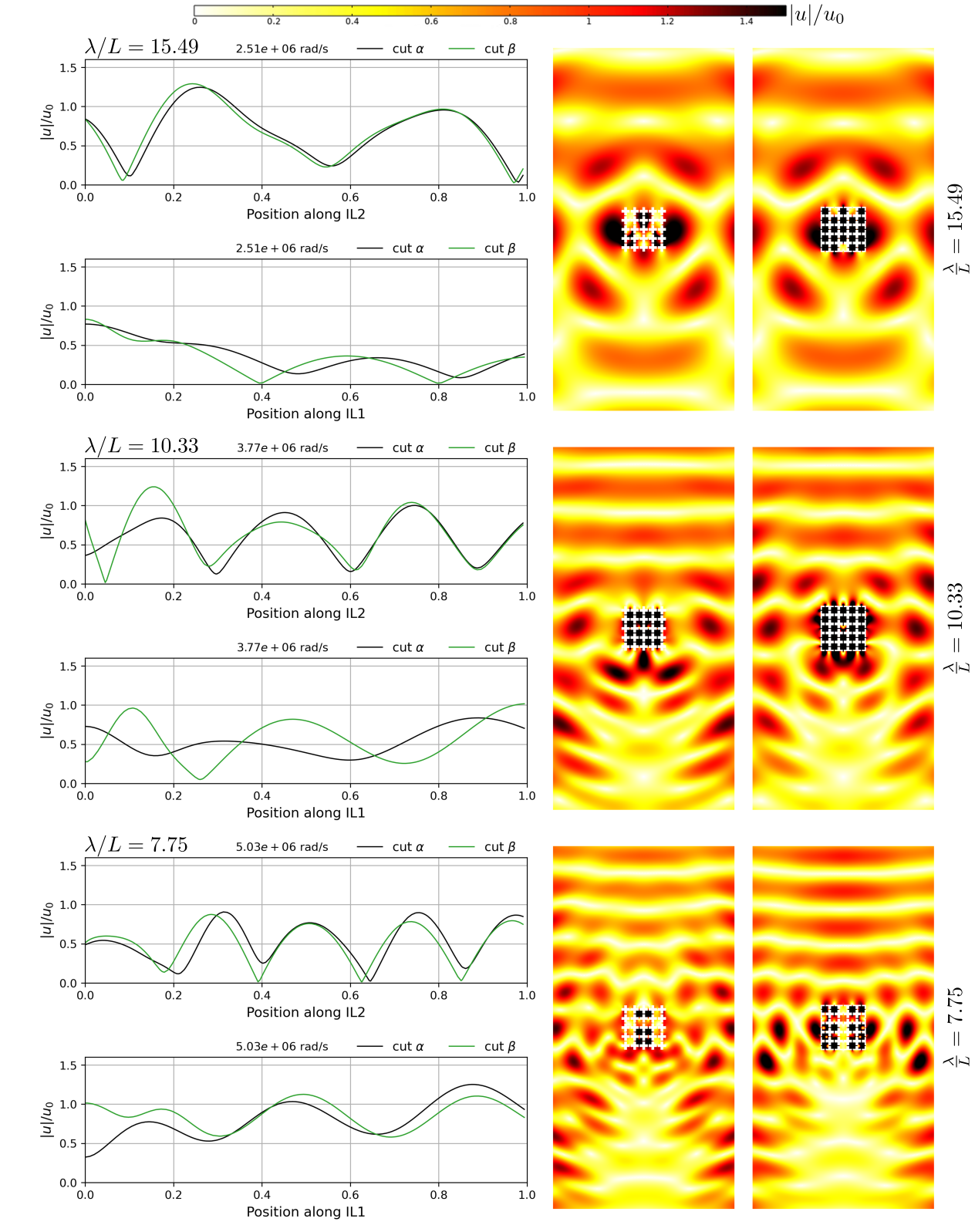}
		\end{subfigure}
		\caption{
		    Differentiated response of the scattering pattern, driven by the boundary contribution of $\alpha$ and $\beta$ interfaces, as a pressure wave interacts with a $5L\times5L$ metamaterial block.
		    The normalized displacement along inspection lines 1 and 2, as well as the normalized displacement field associated with a metamaterial block with $\alpha$- and $\beta$-type boundaries, are shown for the frequencies:
			\textit{(upper)} $2.51\ \mathrm{M\,rad/s}$,
			\textit{(middle)} $3.77\ \mathrm{M\,rad/s}$,
			and \textit{(lower)} $5.03\ \mathrm{M\,rad/s}$.
			The corresponding wavelengths become comparable to the size of the unit cell: small differences between $\alpha$- and $\beta$-cut start being observable.
			The wavelength $\lambda = 10.33L$ is the one at which the first differences between the two cuts start becoming apparent: we will then suppose that at this threshold value and below, separation of scales is not expected to hold true anymore.
		}
		\label{fig:mstd_differences-5_mid}
	\end{figure}

\begin{figure}[H]
		\begin{subfigure}{\textwidth}
			\centering
			\includegraphics[width=\textwidth]{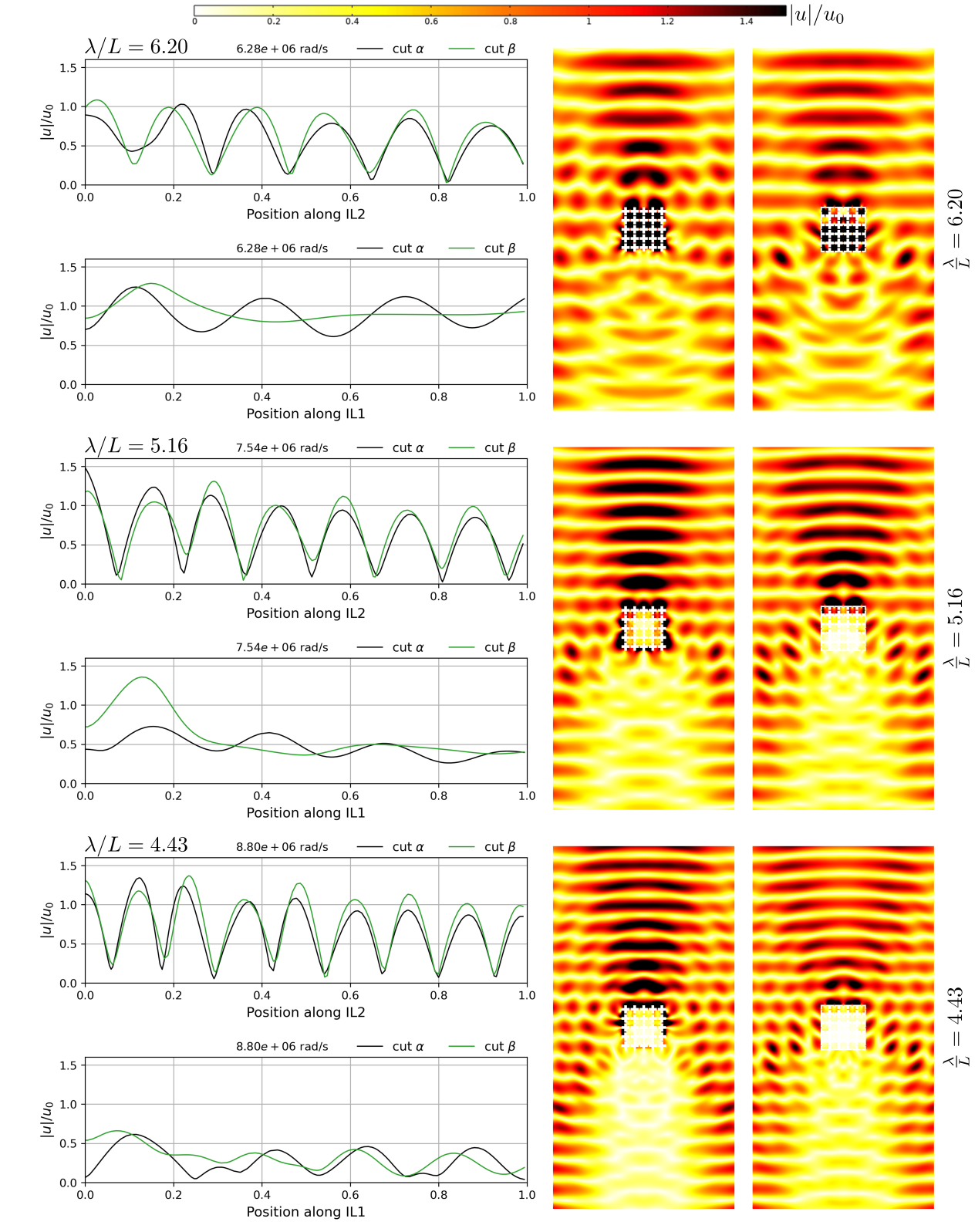}
		\end{subfigure}
		\caption{
		    Differentiated response of the scattering pattern, driven by the boundary contribution of $\alpha$ and $\beta$ interfaces, as a pressure wave interacts with a $5L\times5L$ metamaterial block.
		    The normalized displacement along inspection lines 1 and 2, as well as the normalized displacement field associated with a metamaterial block with $\alpha$- and $\beta$-type boundaries, are shown for the frequencies:
			\textit{(upper)} $6.28\ \mathrm{M\,rad/s}$,
			\textit{(middle)} $7.54\ \mathrm{M\,rad/s}$,
			and \textit{(lower)} $8.80\ \mathrm{M\,rad/s}$.
			The corresponding wavelengths become closer to the unit cell's size: larger differences between $\alpha$- and $\beta$-cuts can be observed. Separation of scales hypothesis breaks down.
		}
		\label{fig:mstd_differences-5_high}
	\end{figure}

\begin{figure}[H]
		\begin{subfigure}{\textwidth}
			\centering
			\includegraphics[width=\textwidth]{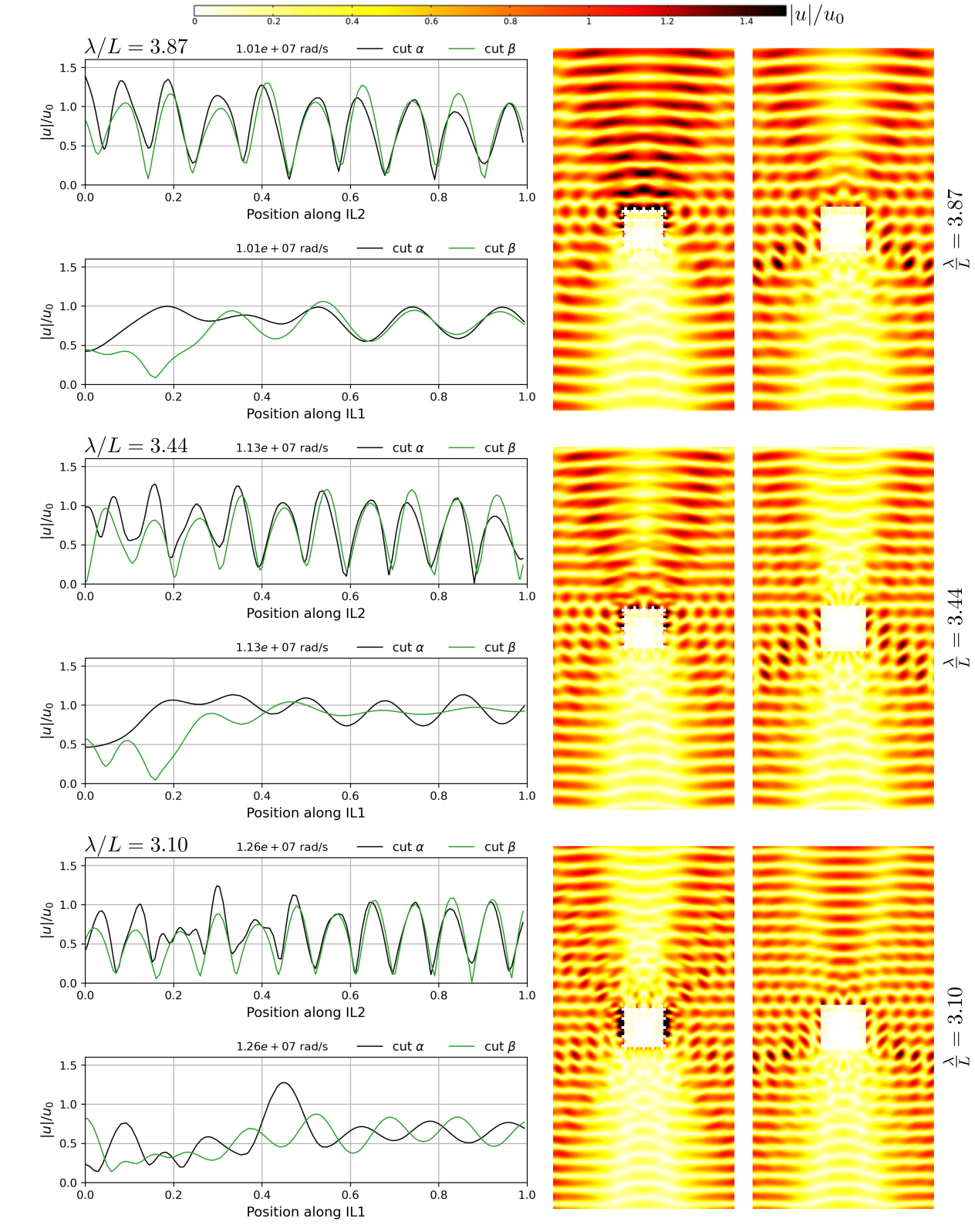}
		\end{subfigure}
		\caption{
		    Differentiated response of the scattering pattern, driven by the boundary contribution of $\alpha$ and $\beta$ interfaces, as a pressure wave interacts with a $5L\times5L$ metamaterial block.
		    The normalized displacement along inspection lines 1 and 2, as well as the normalized displacement field associated with a metamaterial block with $\alpha$- and $\beta$-type boundaries, are shown for the frequencies:
			\textit{(upper)} $10.05\ \mathrm{M\,rad/s}$,
			\textit{(middle)} $11.31\ \mathrm{M\,rad/s}$,
			and \textit{(lower)} $12.56\ \mathrm{M\,rad/s}$.
			The corresponding wavelengths are very close to the unit cell's size: the differences between $\alpha$- and $\beta$-cuts are important. Separation of scales hypothesis does not hold true anymore.
		}
		\label{fig:mstd_differences-5_band_gap}
	\end{figure}

    Figure \ref{fig:geometry_problem} shows the geometry of the considered scattering problem.
    A plane elastic wave impinges on a finite metamaterial block (a square specimen whose side consists of 5 unit cellas) embedded in a homogeneous medium, generating a scattered wave field due to the interaction with the microstructured region.
    The figure also indicates the two inspection directions along which the displacement is evaluated for a quantitative comparison.
    The displacement profiles along these directions are reported in the left panels of Figs. \ref{fig:mstd_differences-5_low}--\ref{fig:mstd_differences-5_band_gap}.
    The right panels of Figs. \ref{fig:mstd_differences-5_low}--\ref{fig:mstd_differences-5_band_gap} instead show the corresponding fully resolved displacement fields obtained from the simulations for the two considered cuts of the metamaterial.
    
    It can be noticed, by direct inspection of Figs. \ref{fig:mstd_differences-5_low}--\ref{fig:mstd_differences-5_band_gap} that when the wavelength is much larger than the characteristic dimensions of the macroscopic specimen, the displacement fields obtained for the two cuts are nearly indistinguishable.
    
    As the frequency increases, however, small differences between the $\alpha$- and $\beta$-cuts begin to appear.
    Small differences start arising already for $\lambda=10.33L$
    
    However, the most consistent differences arise when the wavelength becomes smaller than this threshold value.
    We can then infer that the separation of scale hypothesis starts breaking down when the wavelength $\lambda$ is about 10 times the size of single unit cell.
    In this regime ($\lambda < 10.33L$), the incident wave starts to interact directly with the microscopic geometric details of the interface, leading to noticeable differences in the scattering fields of the two cuts.
    
    These observations clearly show that finite-size metamaterials cannot always be consistently treated within a classical homogenization framework at higher frequencies.
    When the wavelength approaches the characteristic length scale of the system---namely the unit cell size---the dynamic response becomes sensitive to the specific configuration of the interface.
    As a result, specimens that share identical bulk properties but differ only in the way the microstructure is truncated can exhibit different macroscopic responses, highlighting the need for homogenized models capable of accounting for boundary-related effects.
    
    These differences can therefore be interpreted as \textbf{dynamic size effects} arising in finite metamaterial specimens.
    They can be associated with local resonant phenomena occurring near the interface at the level of the underlying microstructure.
    As the wavelength approaches the characteristic size of the unit cell, the boundary region no longer behaves as a homogenized medium, and localized dynamic interactions between the incident wave and the microstructural elements close to the interface become significant.
    
 A similar story for a $20\times20$ macroscopic specimen is given in Appendix: a similar trend can be observed in Figs. \ref{fig:mstd_differences-20_low}--\ref{fig:mstd_differences-20_band_gap}.
    The threshold value $\lambda=10.33L$ for which differences between $\alpha$- and $\beta$-cut start arising is confirmed also in this case.
    In the following, we will then assume that the separation of scale hypothesis starts breaking down when $\lambda=10.33L$.
    
	\newpage

\FloatBarrier

\section{Relaxed Micromorphic Modeling of Finite-Size Metamaterials}
	
	The metamaterials' dispersive behavior is impossible to capture using a homogenized linear-elastic Cauchy model due to its strictly linear dispersion relations \cite{Liu.2018}.
	Conversely, simulations that include all the geometrical details \cite{Wei.2021,dAgostino.2020,AlBabaa.2018,Veres.2013} can precisely simulate the metamaterial's dispersive response using the classical elasticity model, but at the expense of high computational costs \cite{Xue.2023}.
	Therefore, homogenized models \cite{Faraci.2026,Li.2026,Yang.2024,Yang.2021,Boisse.2018} are widely used as a cost effective alternative.
	
	The \rmm has been shown to capture metamaterials' response at frequencies below and in the band gap.
	For instance, Madeo et al. have conducted studies on infinite and finite-size domains \cite{Rizzi.2022b,Madeo.2018,Rizzi.2021,Aivaliotis.2020} that show agreement with the geometrically detailed simulations.
	Moreover, the micromorphic model exhibits a wider frequency range of application in comparison to alternatives such as micropolar and second gradient models \cite{Rizzi.2024}.
	
	The differences recently observed by some of the authors in the traction profiles at the metamaterial's interface in \cite{PerezRamirez.2024} suggest that additional actions are needed to account for differentiated homogenized boundary effects, which evidently occur when considering finite-size metamaterials' specimens with different interfaces. 
	In spite of the uniqueness of the balance equations derived from the first variation of a relaxed micromorphic continuum, we propose that changes in the inertia distribution along the boundary signal the existence of a supplementary microstructure-related interface action.
	We thus propose that a kinetic surface energy can account for the distinct responses of different interfaces.
	The associated modified equilibrium conditions are presented in the following sections.

	\subsection{Reduced Relaxed Micromorphic Model}
	\label{subsec:rrmm}
	
	We pose that in a time interval $[0,\tau]$ and for a body $\Omega$, whose boundaries are $\partial \Omega$, the equilibrium conditions of such body can be found when minimizing its action functional $\mathcal{A}$:
	
	\begin{align}
		\mathcal{A}_{\rm rrm}=\iint\limits_{\Omega \times \left[0,\tau\right]} 
		\mathcal{L} \left(\nabla u, \dot{u}, \nabla \dot{u}, P, \dot{P}\right)
		\, dx \, dt \, .
		\label{eq:act_func_rrmm}
	\end{align}
	
	Here, the Lagrangian reads: 
	
	\begin{align}
		\mathcal{L}_{\rm rrm}
		\coloneqq&
		K \left(\dot{u},\nabla \dot{u}, \dot{P}\right) 
		-
		W \left(\nabla u, P\right)
		\ .
		\label{eq:lag_rrmm} 
	\end{align}
	
	The kinetic $K$ and strain $W$ energies of the \rrmm are:
	
	\begin{align}
		K \left(\dot{u},\nabla \dot{u}, \dot{P}\right) 
		=&
		\dfrac{1}{2}\rho \, \langle \dot{u},\dot{u} \rangle + 
		\dfrac{1}{2} \langle \mathbb{J}_{\rm m}  \, \text{sym} \, \dot{P}, \text{sym} \, \dot{P} \rangle 
		+ \dfrac{1}{2} \langle \mathbb{J}_{\rm c} \, \text{skew} \, \dot{P}, \text{skew} \, \dot{P} \rangle
		\notag
		\\
		&
		+ \dfrac{1}{2} \langle \mathbb{T}_{\rm e} \, \text{sym}\nabla \dot{u}, \text{sym}\nabla \dot{u} \rangle
		+ \dfrac{1}{2} \langle \mathbb{T}_{\rm c} \, \text{skew}\nabla \dot{u}, \text{skew}\nabla \dot{u} \rangle
		\label{eq:kin_rrmm}
		\ ,
		\\
		W \left(\nabla u, P\right)
		=& 
		\dfrac{1}{2} \langle \mathbb{C}_{\rm e} \, \text{sym}\left(\nabla u -  \, P \right), \text{sym}\left(\nabla u -  \, P \right) \rangle
		\notag
		\\
		&
		+ \dfrac{1}{2} \langle \mathbb{C}_{\rm c} \, \text{skew}\left(\nabla u -  \, P \right), \text{skew}\left(\nabla u -  \, P \right) \rangle
		\notag
		\\
		&
		+ \dfrac{1}{2} \langle \mathbb{C}_{\rm micro} \, \text{sym}  \, P,\text{sym}  \, P \rangle
		\label{eq:str_rrmm} 
		\ ,
	\end{align}
	
	\noindent
	where $u \in \mathbb{R}^{3}$ is the macroscopic displacement field, $P \in \mathbb{R}^{3\times 3}$ is the non-symmetric micro-distortion tensor, $\rho$ is the macroscopic apparent density, $\mathbb{J}_{\rm m}$, $\mathbb{J}_{\rm c}$, $\mathbb{T}_{\rm e}$, $\mathbb{T}_{\rm c}$ are 4th order micro-inertia tensors, and $\mathbb{C}_{\rm e}$, $\mathbb{C}_{\rm m}$, $\mathbb{C}_{\rm c}$ are 4th order elasticity tensors.
	The elasticity, micro-distorsion and micro-inertia tensors in Voigt notation are given by:
	
	\begin{align}
		\mathbb{C}_{\rm e}
		&= 
		\begin{pmatrix}
			\kappa_{\rm e} + \mu_{\rm e}	& \kappa_{\rm e} - \mu_{\rm e}				& \star & \dots	& 0\\ 
			\kappa_{\rm e} - \mu_{\rm e}	& \kappa_{\rm e} + \mu_{\rm e} & \star & \dots & 0\\
			\star & \star & \star & \dots & 0\\
			\vdots & \vdots	& \vdots & \ddots &\\ 
			0 & 0 & 0 & & \mu_{\rm e}^{*}
		\end{pmatrix},
		&
		\mathbb{C}_{\rm micro}
		&=
		\begin{pmatrix}
			\kappa_{\rm m} + \mu_{\rm m}	& \kappa_{\rm m} - \mu_{\rm m}				& \star & \dots	& 0\\ 
			\kappa_{\rm m} - \mu_{\rm m}	& \kappa_{\rm m} + \mu_{\rm m} & \star & \dots & 0\\
			\star & \star & \star & \dots & 0\\
			\vdots & \vdots	& \vdots & \ddots &\\ 
			0 & 0 & 0 & & \mu_{\rm m}^{*}
		\end{pmatrix},
		\notag
		\\[2mm]
		\mathbb{J}_{\rm m}
		&=
		\rho L_{\rm c}^2
		\begin{pmatrix}
			\kappa_\gamma + \gamma_{1} & \kappa_\gamma - \gamma_{1} & \star & \dots & 0\\ 
			\kappa_\gamma - \gamma_{1} & \kappa_\gamma + \gamma_{1} & \star & \dots & 0\\ 
			\star & \star & \star & \dots & 0\\
			\vdots & \vdots & \vdots & \ddots &\\ 
			0 & 0 & 0 & & \gamma^{*}_{1}\\ 
		\end{pmatrix},
		&
		\mathbb{T}_{\rm e}
		&=
		\rho L_{\rm c}^2
		\begin{pmatrix}
			\overline{\kappa}_{\gamma} + \overline{\gamma}_{1} & \overline{\kappa}_{\gamma} - \overline{\gamma}_{1} & \star & \dots	& 0\\ 
			\overline{\kappa}_{\gamma} - \overline{\gamma}_{1} &   \overline{\kappa}_{\gamma} + \overline{\gamma}_{1} & \star & \dots & 0\\ 
			\star & \star & \star & \dots & 0\\
			\vdots & \vdots & \vdots & \ddots &\\ 
			0 & 0 & 0 & & \overline{\gamma}^{*}_{1}
		\end{pmatrix},
		\label{eq:tensors}
	\end{align}
	
	\begin{align}
		\mathbb{J}_{\rm c}
		&=
		\rho L_{\rm c}^2
		\begin{pmatrix}
			\star & 0 & 0\\ 
			0 & \star & 0\\ 
			0 & 0 & 4\,\gamma_{2}
		\end{pmatrix},
		&
		\mathbb{T}_{\rm c}
		&=
		\rho L_{\rm c}^2
		\begin{pmatrix}
			\star & 0 & 0\\ 
			0 & \star & 0\\ 
			0 & 0 & 4\,\overline{\gamma}_{2}
		\end{pmatrix},
		&
		\mathbb{C}_{\rm c}
		&= 
		\begin{pmatrix}
			\star & 0 & 0\\ 
			0 & \star & 0\\ 
			0 & 0 & 4\,\mu_{\rm c}
		\end{pmatrix}.
		\notag
	\end{align}
	
	For the case in which there are $f$ external boundary forces and no body forces , the first variation of the action functional $\mathcal{A}$ leads to the following strong form of equilibrium equations:
	
	\begin{align}
		\overline{\sigma}
		&= 
		\widetilde{\sigma}
		-
		s
		&
		\quad
		&\text{in}\ \Omega
		\notag
		\\
		\text{Div}\widetilde{\sigma}
		&=
		\rho\,\ddot{u}
		-
		\text{Div}\widehat{\sigma}
		&
		\quad
		&\text{in}\ \Omega
		\notag
		\\
		t
		=
		f
		\quad
		&
		\mathrm{and}
		\quad
		u
		=
		\overline{u}
		&
		\quad
		&\text{on}\ \partial \Omega
		\ ,
		\notag
	\end{align}
	
	\noindent
	where:
	
	\begin{align}
		\widetilde{\sigma}
		&\coloneqq
		\mathbb{C}_{\rm e}\,\text{sym}(\nabla u-P) + \mathbb{C}_{\rm c}\,\text{skew}(\nabla u-P)
		\,,
		&
		\widehat{\sigma}
		&\coloneqq
		\mathbb{T}_{\rm e}\,\text{sym}\nabla\ddot{u} + \mathbb{T}_{\rm c}\,\text{skew}\nabla\ddot{u}
		\,,
		\notag
		\\
		\overline{\sigma}
		&\coloneqq
		\mathbb{J}_{\rm m}\,\text{sym}\ddot{P} + \mathbb{J}_{\rm c}\,\text{skew}\ddot{P}
		\,,
		&
		s
		&\coloneqq
		\mathbb{C}_{\rm micro}\, \text{sym} P
		\,,
		\notag
		\\
		t
		&\coloneqq
		\left( \widetilde{\sigma} + \widehat{\sigma} \right) \, n 
		\,.
		&
		&
		\notag
	\end{align}

	\subsection{Towards Breaking Separation of Scales in homogenized models: New Interface Inertia Term}
	
	Simulations and experiments of finite-size metamaterials show that cutting the metamaterial's specimen at different locations changes the response of the system for sufficiently high frequencies \cite{Hermann.2026}.	
	In fact, the mass distribution close to the interface is different according to where the underlying structure is cut (\fig{fig:interface_mass_distribution}).
	We propose that the differentiated response in the band gap is driven by the inertial contribution of the region close to the boundary.
	Therefore, to capture the boundary effects, we enhance the micromorphic model using a \textit{kinetic surface energy}:\footnotemark
	
	\footnotetext{We introduce the "-" sign in the definition of the kinetic energy to keep the convention on the sign for the forces acting at the interfaces: we want here the same sign as the externally applied forces "+" (see eq. \ref{eq:bc_total}).}
	
	\begin{equation}
		K_{\partial \Omega} \left(\dot{u}\right) 
		\coloneqq
		-\dfrac{1}{2}\rho_{\partial \Omega} \, \langle \dot{u},\dot{u} \rangle
		\,,
	\end{equation}
	
	\noindent
	where $\rho_{\partial \Omega}$ depends on the choice of truncating plane.
	
	This leads to the total action functional:
	
	\begin{align}
		\mathcal{A}_{\rm total}
		&=
		\mathcal{A}_{\partial \Omega}
		+
		\mathcal{A}_{\rm rrm}
		\notag
		\\
		&=
		\iint\limits_{\partial \Omega \times \left[0,\tau\right]} 
		K_{\partial \Omega} \left(\dot{u}\right)
		\, dx \, dt \,
		+
		\iint\limits_{\Omega \times \left[0,\tau\right]} 
		\mathcal{L}_{\rm rrm} \left(\nabla u, \dot{u}, \nabla \dot{u}, P, \dot{P}\right)
		\, dx \, dt \, .
		\label{eq:act_func_tot}
	\end{align}

	\subsubsection{Variation of the Action Functional and Equations in Strong Form}
	
	The first variation of $\mathcal{A}_{\partial \Omega}$ after integration by parts is:
	
	\begin{align}
		\partial \mathcal{A}_{\partial \Omega}
		=
		\iint\limits_{\partial \Omega \times \left[0,\tau\right]} 
		\dyniprod{\rho_{\partial \Omega}\,\ddot{u} \,,\, \delta u}
		\, dx \, dt \,
	\end{align}
	
	We then have that the boundary terms of the first variation of the total action functional are:
	
	\begin{align}
		\partial \mathcal{A}_{\rm total}
		=
		\iint\limits_{\partial \Omega \times \left[0,\tau\right]} 
		\dyniprod{\rho_{\partial \Omega}\,\ddot{u} \,,\, \delta u}
		\, dx \, dt \,
		-
		\iint\limits_{\partial \Omega \times \left[0,\tau\right]} 
		\dyniprod{\left( \widetilde{\sigma} + \widehat{\sigma} \right) \, n  \,,\, \delta u}
		\, dx \, dt \,
	\end{align}
	
	As a result, $K_{\partial \Omega}$ adds a boundary term to the strong form of equilibrium equations:
	
	\begin{align}
		\overline{\sigma}
		&= 
		\widetilde{\sigma}
		-
		s
		&
		\quad
		&\text{in}\ \Omega
		\\
		\text{Div}\widetilde{\sigma}
		&=
		\rho\,\ddot{u}
		-
		\text{Div}\widehat{\sigma}
		&
		\quad
		&\text{in}\ \Omega
		\\
		t
		&=
		f
		+
		\rho_{\partial \Omega}\,\ddot{u}
		&
		\quad
		&\text{on}\ \partial \Omega
		\label{eq:bc_total}
		\ .
	\end{align}
	
	The result of \equ{eq:bc_total} introduces an energy-based non-coherent boundary condition of the elastic-interface type discussed in \cite{PerezRamirez.2024,Murdoch.1976,Moeckel.1975,Fried.2007}.

	\subsection{Determination of the Relaxed Micromorphic Parameters}
	\label{subsec:parameters}
	
	We recall here the procedure used to determine the relaxed micromorphic bulk parameters of the "cross-cell" metamaterial given in \fig{fig:unit_cell-characteristics}.
	
	The elasticity and inertia tensors defined in \equ{eq:tensors} take the values given in Table \ref{tab:parameters_RM} when considering the cross-cell defined in \fig{fig:unit_cell-characteristics}. 
	
    \begin{table}[h!]
        \centering
        \begin{tabular}{cccccccc} 
            $\rho$                                  & $\kappa_{\rm m}$              & $\mu_{\rm m}$                     & $\mu_{\rm m}^{*}$                 &
            $\mu_{\rm c}$                           & $\kappa_{\rm e}$              & $\mu_{\rm e}$                     & $\mu_{\rm m}^{*}$                 \\
            kg/m$^3$                                & GPa                           & GPa                               & GPa                               &
            GPa                                     & GPa                           & GPa                               & GPa                               \\
            \hline\hline
            $1485$                                  & $18.91$                       & $7.496$                           & $356.2$                           &
            $0.1000$                                & $12.83$                       & $27.85$                           & $0.6271$                          \\
            \\
            $\overline{\kappa}_{\gamma}$            & $\overline{\gamma}_{1}$       & $\overline{\gamma}_{2}$           & $\overline{\gamma}^{*}_{1}$       &
            $\kappa_\gamma$                         & $\gamma_{1}$                  & $\gamma_{2}$                      & $\gamma^{*}_{1}$                  \\
            \hline\hline
            $0.01186$                               & $0.02536$                     & $0.1752$                          & $0.01980$                         &
            $0.07034$                               & $0.08667$                     & $0.003075$                        & $0.8750$                          \\
        \end{tabular}%
        \caption{
            Parameters of the relaxed micromorphic model for the "cross-cell" metamaterial.
        }
        \label{tab:parameters_RM}
    \end{table}

	\subsubsection{Determination of the Bulk Parameters}
	
The parameters of the relaxed micromorphic model are identified using the two-level optimization procedure introduced in \cite{Sarhil.2026}. First, the static parameters are identified by fitting the static size effects \cite{Sarhil.2024}. Subsequently, the dynamic parameters are calibrated by fitting the dispersion curves of the relaxed micromorphic model to those obtained from the standard Bloch–Floquet analysis.
			
	The relaxed micromorphic model behaves in the static case as a two-scale elasticity model with upper and lower bounds corresponding to standard linear elasticity with elasticity tensors $\Cmicro$ and $\Cmacro$ \cite{Sarhil.2023}. The static parameters of the RMM are identified by fitting the energy of different loading cases of the RMM with the full-fidelity microstructured solution.  To this end, a cost function $r^2_\textrm{static}$ is minimized over several loading cases, indexed by $i$, and for domains containing different numbers $n$ of unit cells as  
	
	\begin{equation}
r^2_\textrm{static} =  \sum_{n}  \sum_{i} \lvert   \int\limits_{\Omega} 
	W^\textrm{het}_{\{i , n\}} (\nabla u) \, dx 	- \int\limits_{\Omega}  W \left(\nabla u, P, \textrm{Curl} P \right)
		\, dx  \,  \rvert^2 \,  \rightarrow \textrm{min} ,     
\end{equation} 

where $W$ is strain energy density of the homogeneous relaxed micromorphic domain\footnote{The strain energy $W$ for the RMM is given in \equ{eq:str_rrmm}, the energy considered for the static fitting also includes $\textrm{Curl} P$ term in order to account for static size effects. Since those effects become negligible in the dynamic regime, we neglect the $\textrm{Curl} P$ term in the reminder of this paper.} and $W^\textrm{het}$ is the strain energy of the heterogeneous linear-elastic domains.  
For the static identification, the curvature term incorporating $\textrm{Curl}P$ is activated to recover an upper bound, characterized  by linear elasticity  with elasticity tensor $\mathbb{C}_{\mathrm{micro}}$. This requires enforcing the so-called consistent boundary condition to ensure that the intended upper bound is reproduced \cite{Sarhil_consist_2022,Sarhil_consist_2023}. A gradient-based optimization with a line-search algorithm is introduced for second-order loading scenarios in \cite{Sarhil.2024} and first-order loading scenarios in \cite{Sarhil.2026} which we also adopt in this work. For the numerical aspects related to the  $H(\textrm{Curl})$-conforming finite element formulation, we refer to \cite{Sarhil.2021,Schroeder.2022}.

\begin{figure}[!ht]
		\begin{subfigure}{\textwidth}
			\centering
			\includegraphics[width=\textwidth]{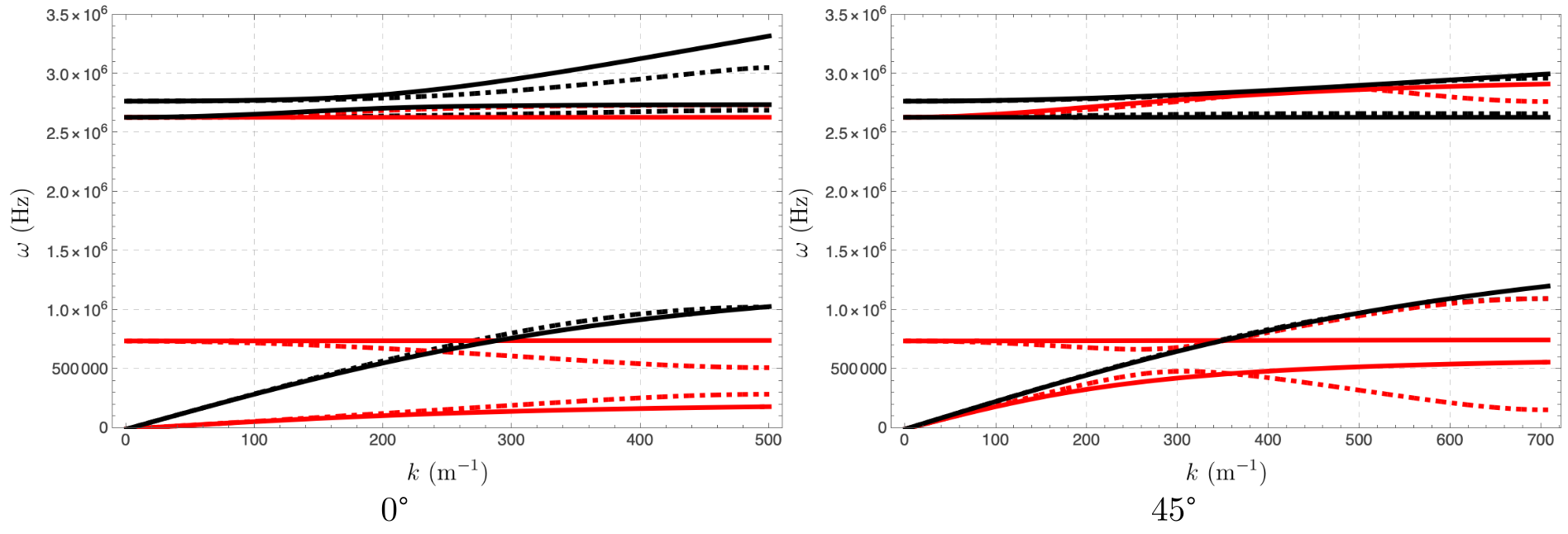}
		\end{subfigure}
		\caption{
		    Dispersion curves of the "cross-cell" metamaterial.
		    Dashed lines show the Bloch-Floquet solution of the microstructured simulation, and solid lines show the corresponding relaxed micromorphic model curves using the parameters of Table \ref{tab:parameters_RM}.
		    Fitting was performed using the incidence angles $0\degree$, and $45\degree$. 
		    Pressure and shear wave curves are shown in black and red, respectively.
		}
		\label{fig:disp_curves}
	\end{figure}

	The dynamic parameters are calibrated by minimizing the discrepancy between the dispersion curves obtained analytically for the relaxed micromorphic model and those computed using Bloch-Floquet simulations on a detailed heterogeneous unit cell.  To this end, an optimization procedure is formulated by minimizing  a cost function $r^2_\textrm{dynamic} $ that fits the first six dispersion branches for two propagation directions, corresponding to the incident angles $\alpha = 0^\circ$ and $\alpha = 45^\circ$, such that

	\begin{equation}
r^2_\textrm{dynamic} =  \sum_{\alpha=0^\circ,45^\circ} \sum_{i=1}^6  \int_k [ r_i (\omega_i^{\alpha} (k)-\omega_i^{\alpha,\textrm{het}}(k))]^2 \textrm{d}k \,  \rightarrow \textrm{min} \,, 
\end{equation}

where $r_i$ are weighting factors which are set higher for the dispersion curves under than band-gap. The results of the fitting procedure are shown in \fig{fig:disp_curves}.

	\subsubsection{Calibration of the Interface Surface Density}
	\label{sec:calibration}
	
	In this subsection we calibrate the interface density parameter $\rho_{\partial \Omega}$ for the $\alpha$- and $\beta$-cuts by using the benchmark problem defined in section \ref{sec:influence_boundary_truncation}.
	
\begin{figure}[H]
		\begin{subfigure}{\textwidth}
			\centering
			\includegraphics[width=0.9 \textwidth]{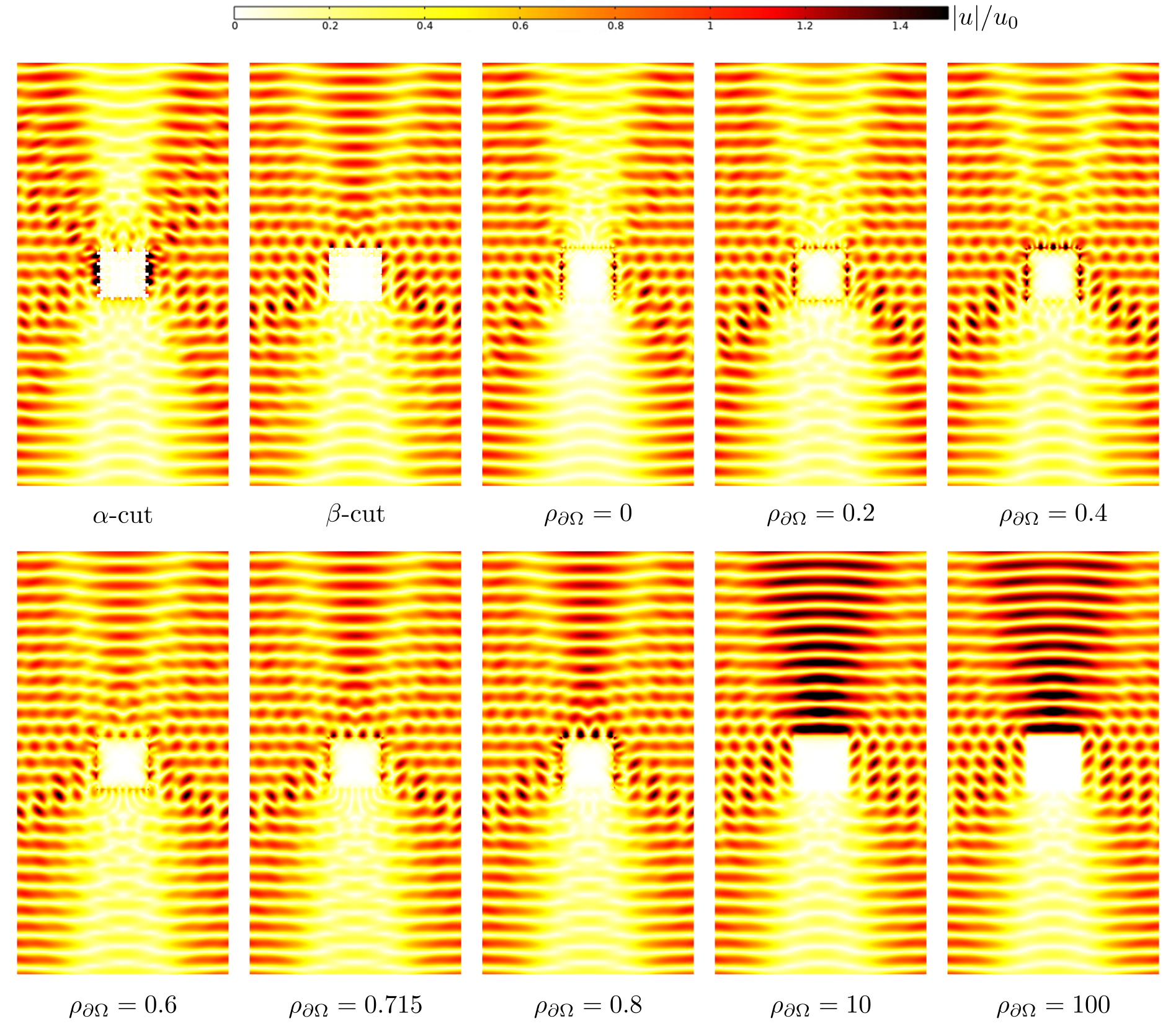}
		\end{subfigure}
		\caption{
		    Calibration of the surface inertia for the $\beta$-cut using a $5L\times5L$ metamaterial block as it interacts with a pressure plane wave that propagates vertically downwards at $12.56$ M\,rad/s and ($\lambda = 3.10L$).
		    From left to right, the normalized displacement field corresponds to the scattering pattern of the:
		    \textit{$\alpha$-cut} microstructured simulation,
		    \textit{$\beta$-cut} microstructured simulation,
		    micromorphic simulations using the surface inertia contribution and the indicated surface density value.
			By direct observation, we can infer that the case $\rho_{\partial \Omega}=0.715$ $\rm kg/m^2$ retrieves the $\beta$-cut results and $\rho_{\partial \Omega}=100$ $\rm kg/m^2$ shows the limiting case where the interface does not move.
		}
		\label{fig:5_calibration}
	\end{figure}

	The choice of the $5\times5$ specimen is motivated by the fact that boundary eﬀects are most pronounced in smaller metamaterial domains.
	Since the objective of the calibration is to identify a surface density that captures the dynamic influence of the interface, it is advantageous to consider a configuration in which the response is strongly aﬀected by the boundary.
	In larger specimens the relative influence of the interface is weaker, which would make the identification of the surface parameter less effective.
	We also choose the calibration frequency $\omega=12.56$ M\,rad/s because it corresponds to a regime in which the wavelength of the incident wave is comparable to only a few unit cells.
	As shown in the previous subsection, this is precisely the regime where the separation of scales hypothesis breaks down and the differences between the $\alpha$- and $\beta$-cuts become clearly visible.
	Calibrating the surface density in this frequency range therefore allows the additional boundary inertia introduced in the relaxed micromorphic model to capture the interface-driven dynamic eﬀects observed in the fully resolved simulations.
	Figure \ref{fig:5_calibration} shows the displacement fields obtained from the fully resolved simulations for the $\alpha$- and $\beta$-cuts together with the relaxed micromorphic predictions for different values of the surface density $\rho_{\partial \Omega}$.
	
	By direct inspection of Fig. \ref{fig:5_calibration} we can infer that the value $\rho_{\partial \Omega}=0$ can be assumed to correspond to the $\alpha$-cut specimen.
	
	The interface density for the $\alpha$-cut is taken to be zero because the standard relaxed micromorphic model already reproduces the response of the $\alpha$-truncated specimen without any additional	boundary correction.
	In other words, $\evalat{\surfdens}{\alpha}=0$ does not mean that the physical $\alpha$-interface has no mass.
	Rather, it means that the inertial effect of the $\alpha$-boundary is already effectively contained in the homogenized description with coherent boundary conditions.
	Therefore, no extra surface inertia is needed to match the fully resolved $\alpha$-cut simulation.
	
	The surface density $\surfdens$ should be interpreted as a relative interface correction with respect to the reference homogenized boundary response.
	In the present calibration, the $\alpha$-cut plays the role of this reference configuration.
	The relaxed micromorphic model with $\surfdens=0$ naturally matches the $\alpha$-cut response.
	By contrast, the $\beta$-cut generates a different boundary mass distribution, which is not captured by the standard coherent-interface formulation.
	It therefore requires an additional effective surface density, namely $\evalat{\surfdens}{\alpha}=0.715$ $\rm kg/m^2$.
	
	Thus, the meaning of $\evalat{\surfdens}{\alpha}=0$ is that the $\alpha$-cut does not introduce an additional truncation-induced inertial correction beyond the one already represented by the bulk relaxed micromorphic model and its standard boundary condition.
	The non-zero value assigned to the $\beta$-cut measures how much its boundary dynamics deviates from this reference interface.
	By progressively adjusting the parameter $\rho_{\partial \Omega}$, we then find the best agreement with the fully resolved $\beta$-cut simulation to be $\rho_{\partial \Omega}=0.715$ $\rm kg/m^2$ (see Fig. \ref{fig:5_calibration} ).
	
	For completeness, Fig. \ref{fig:5_calibration} also shows the limiting case corresponding to very large negative values of the surface density.
	In this limit the inertial contribution of the interface dominates the boundary condition, effectively preventing motion of the interface and leading to a nearly rigid boundary response.
	This limiting behavior is already almost reached for $\rho_{\partial \Omega} = 10$, while the case $\rho_{\partial \Omega} = 100$ essentially represents the asymptotic rigid-interface limit.

\FloatBarrier

\section{Results}
	
	In this section we show that the calibrated value for the interface density $\surfdens=0.715 \rm kg/m^2$ is able to adjust the relaxed micromorphic response for a large frequency range and not only for the frequency used in the calibration procedure described in section \ref{sec:calibration}.
	This means that our hypothesis according to which the metamaterial's response is strongly affected by lower-scale vibration mechanisms across the interface is fully justified.

	\subsection{Frequency- and Size-Validation of the Interface Density Calibrated for the $\beta$-Cut}
	
	We now assess the predictive capability of the interface density calibrated on the $\beta$-cut when it is used at frequencies different from the calibration frequency.
	The calibrated value of the interface density is $\surfdens=0.715 \rm kg/m^2$.
	This value was identified at a frequency for which the wavelength is comparable to only a few unit cells and, therefore, the classical separation of scales hypothesis does not hold.
	The objective of the present validation is to verify whether the same value of $\surfdens$ is able to reproduce the scattering response of the $\beta$-Cut over a broader frequency range and for different specimen sizes.
	
	The comparison is performed for three finite metamaterial blocks of size $5L\times5L$, $10L\times10L$L, and $20L\times20L$.
	For each size, the fully resolved microstructured simulations corresponding to the $\alpha$-cut and $\beta$-cut are compared with relaxed micromorphic simulations obtained either with $\surfdens=0$, corresponding to the standard coherent-interface formulation, or with $\surfdens=0.715 \rm kg/m^2$, corresponding to the interface-density correction calibrated on the $\beta$-cut.
	
	At the lowest frequencies, for which the wavelength remains much larger than the characteristic size of the unit cell, the separation of scales hypothesis is still satisfied.
	In this regime, the scattering patterns associated with the $\alpha$-cut and $\beta$-cut are very similar for all three specimen sizes.
	The microstructured simulations therefore do not display a significant sensitivity to the specific boundary truncation, and the relaxed micromorphic model provides an appropriate homogenized response.
	Moreover, the introduction of the interface density has practically no visible effect on the relaxed micromorphic solution: the results obtained with $\surfdens=0$ and with $\surfdens=0.715 \rm kg/m^2$ are nearly indistinguishable.
	This confirms that the surface inertial contribution is naturally inactive in the low-frequency regime, where the boundary details are not resolved by the propagating wave.
	
	When the frequency increases and the wavelength reaches approximately $\lambda=10.33L$, the separation of scales starts to break down, as observed in Section \ref{sec:influence_boundary_truncation}.
	In the frequency interval extending from $\lambda=10.33L$ down to approximately $\lambda=6.20L$, the $\alpha$-cut and $\beta$-cut begin to produce different scattering responses.
	In this intermediate regime, the relaxed micromorphic simulations remain qualitatively consistent with the corresponding fully resolved solutions, but quantitative deviations are still observed.
	Nevertheless, the introduction of the calibrated interface density systematically improves the agreement with the $\beta$-cut response.
	In other words, although the corrected relaxed micromorphic model does not fully remove all discrepancies in this frequency range, the use of $\surfdens=0.715 \rm kg/m^2$ always moves the homogenized response closer to the fully resolved $\beta$-cut solution than the standard model with $\surfdens=0$.
	
	The fact that these quantitative differences persist for the $5L\times5L$, $10L\times10L$L, and $20L\times20L$ specimens suggests that, in this transition regime, additional interface mechanisms may also contribute to the scattering response; their size-independence indicates that this is not a limitation of homogenized modeling itself, but rather of an incomplete representation of the interface effects. In particular, effects such as surface elasticity or anisotropic interface inertia may become relevant, especially because multiple scattering inside the finite metamaterial block can couple pressure and shear waves.
	These additional interface enhancements are, however, beyond the scope of the present work.
	
	For higher frequencies, corresponding to wavelengths $\lambda=5.16L$ and smaller, the influence of the boundary truncation becomes much more pronounced.
	In this regime, the fully resolved $\alpha$-cut and $\beta$-cut simulations exhibit clearly distinct scattering patterns.
	The relaxed micromorphic model with $\surfdens=0$ accurately reproduces the response of the $\alpha$-cut, while the relaxed micromorphic model enhanced with the calibrated value $\surfdens=0.715 \rm kg/m^2$ recovers very well the response of the $\beta$-cut.
	The agreement is very good both qualitatively, in terms of the overall structure of the scattered field, and quantitatively, in terms of the displacement amplitudes and spatial distribution.
	
	These results show that a single scalar value of the interface density, calibrated at one frequency where separation of scales is lost, remains effective over a broad high-frequency range and for different specimen sizes.
	The surface inertia term therefore provides a robust correction to the relaxed micromorphic model, allowing it to distinguish between two finite metamaterial specimens that share the same bulk properties but differ only in the truncation of their external boundaries.
	
\newpage
	
	\begin{table}[htbp]
      \centering
      \begin{threeparttable}
        \caption{Relative $L^2$ error [\%] of the normalized displacement field of the domain around the metamaterial block. The relative error of the micromorphic simulations is reported with respect to the $\alpha$- and $\beta$-cut types of boundaries of the corresponding microstructured simulation, at the prescribed angular frequencies for $N=5$.}
        \label{tab:frequency-errors-N5}
        \footnotesize
        \setlength{\tabcolsep}{2.7pt}
        \renewcommand{\arraystretch}{1.14}
        \begin{tabular}{@{}l S[table-format=-1.3]
          *{10}{S[table-format=2.2]} S[table-format=2.2]@{}}
          \freqheader
          $\alpha$ &  0     & 4.47 &  8.69 & 26.39 & 33.73 & 17.25 & 19.82 & 26.96 & 22.64 & 28.30 & 27.12 & 21.54 \\
          $\alpha$ & 0.715 & 4.64 &  8.29 & 26.94 & 38.04 & 25.26 & 30.29 & 35.91 & 34.25 & 32.65 & 42.45 & 27.87 \\
          \addlinespace[2pt]
          $\beta$  &  0     & 3.96 & 12.14 & 26.85 & 33.81 & 24.41 & 25.81 & 28.29 & 32.12 & 36.39 & 37.84 & 26.16 \\
          $\beta$  & 0.715 & 5.53 & 18.72 & 28.83 & 27.87 & 14.25 &  9.74 &  9.56 & 10.37 & 10.77 & 12.98 & 14.86 \\
          \bottomrule
        \end{tabular}
      \end{threeparttable}
    \end{table}

\newpage
	
\begin{figure}[H]
		\centering
		\begin{subfigure}{\textwidth}
			\centering
			\includegraphics[width=\textwidth]{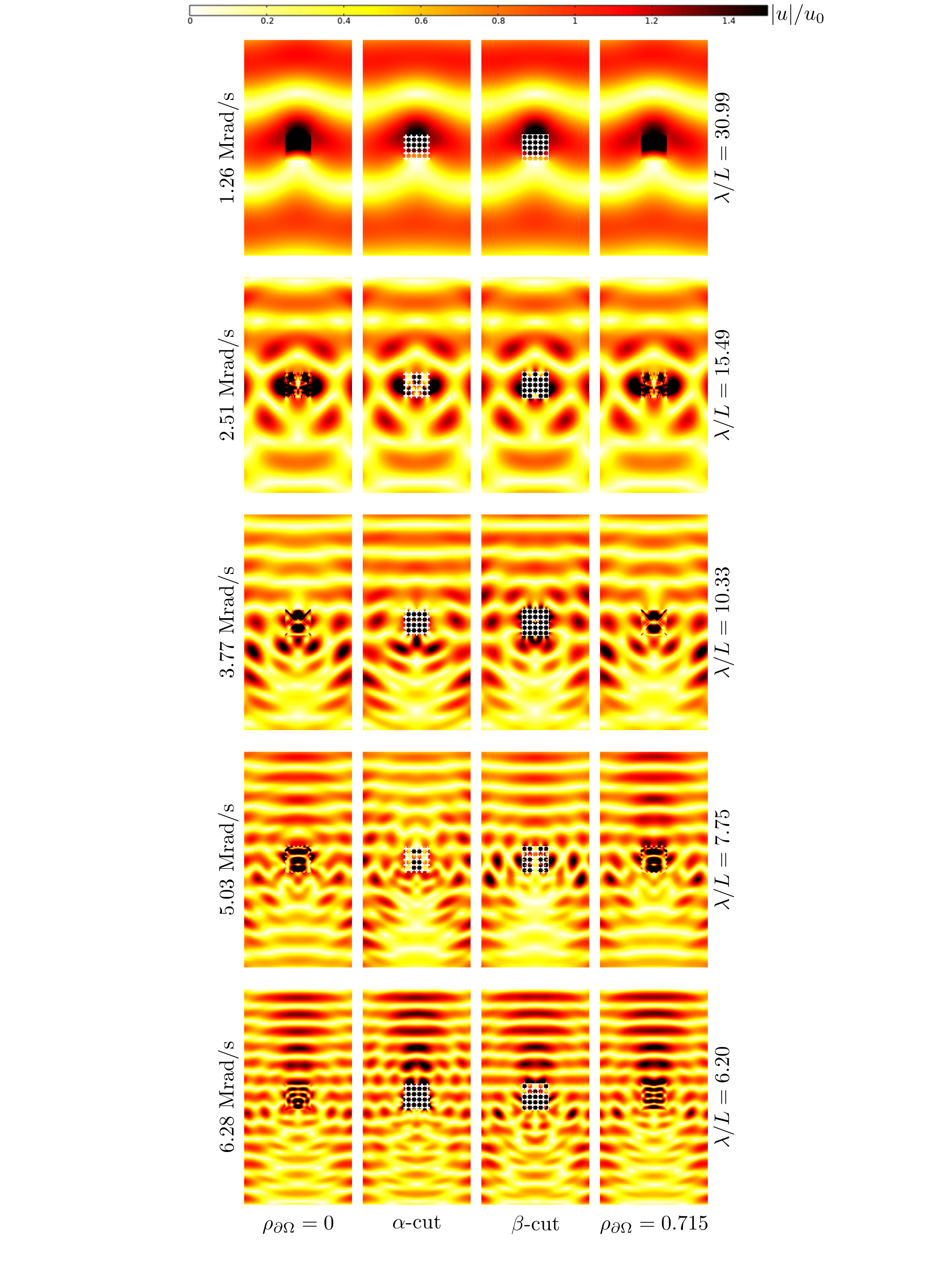}
		\end{subfigure}
		\caption{
		    Response at lower frequencies of a $5L\times5L$ metamaterial block as it interacts with a pressure plane wave that propagates vertically downwards.
		    The normalized displacement fields shown in each column correspond to the resulting scattering pattern for five types of implementation, from left to right: 
		    \textit{($\rho_{\partial \Omega}=0$)} micromorphic simulations without interface contribution,
		    \textit{($\alpha$-cut)} microstructured simulation of the architected material using $\alpha$-type boundaries,
		    \textit{($\beta$-cut)} microstructured simulation of the architected material using $\beta$-type boundaries, and
		    \textit{($\rho_{\partial \Omega}=0.715$)} enhanced micromorphic simulations using the surface inertia contribution calibrated for the $\beta$-cut.
		}
		\label{fig:size_comp_N_5_1}
	\end{figure} 
	
\newpage
	
    \begin{figure}[H]
		\centering
		\begin{subfigure}{\textwidth}
			\centering
			\includegraphics[width=\textwidth]{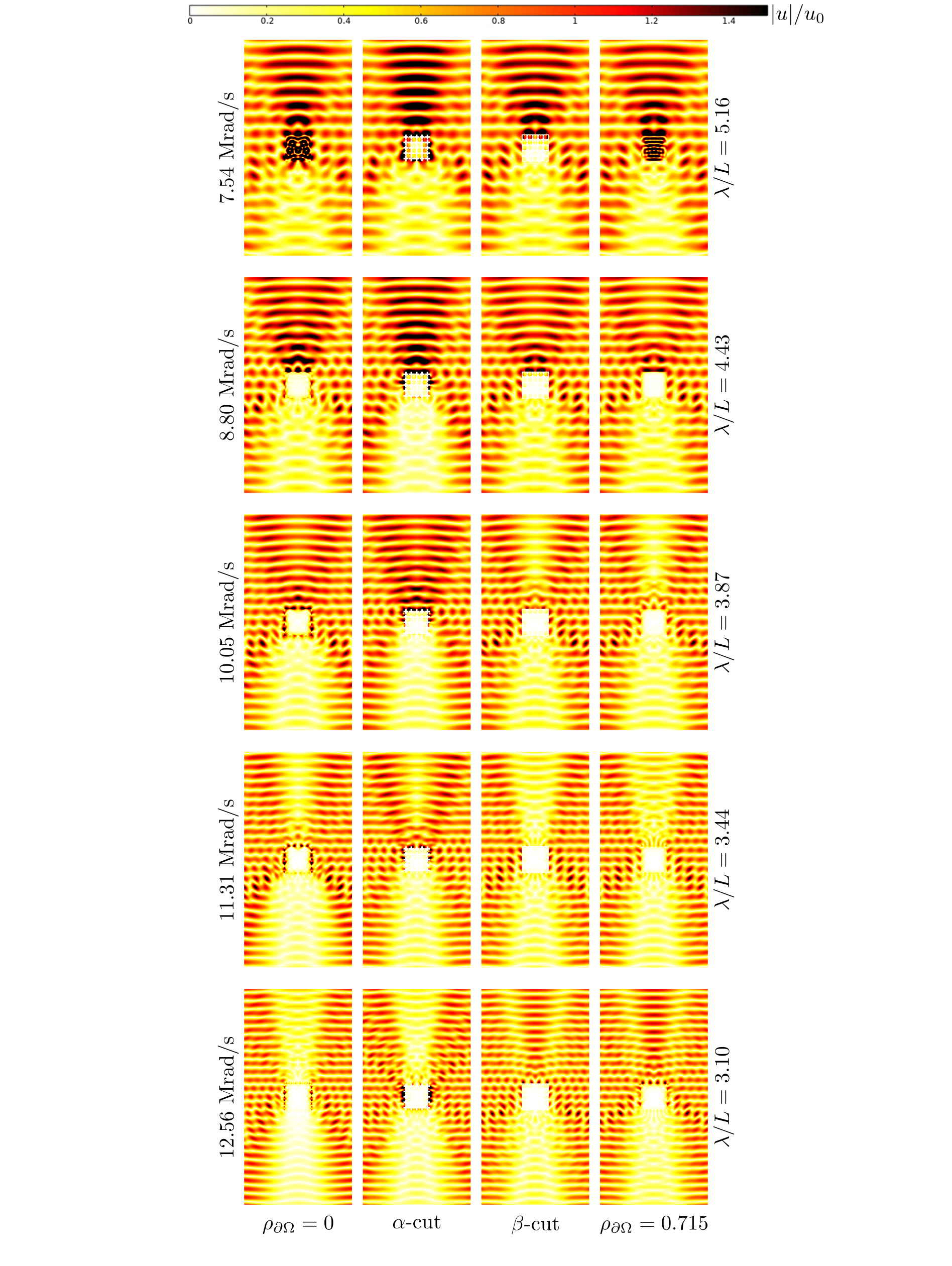}
		\end{subfigure}
		\caption{
		    Response at higher frequencies of a $5L\times5L$ metamaterial block as it interacts with a pressure plane wave that propagates vertically downwards.
		    The normalized displacement fields shown in each column correspond to the resulting scattering pattern for five types of implementation, from left to right:
		    \textit{($\rho_{\partial \Omega}=0$)} micromorphic simulations without interface contribution,
		    \textit{($\alpha$-cut)} microstructured simulation of the architected material using $\alpha$-type boundaries,
		    \textit{($\beta$-cut)} microstructured simulation of the architected material using $\beta$-type boundaries, and
		    \textit{($\rho_{\partial \Omega}=-0.715$)} enhanced micromorphic simulations using the surface inertia contribution calibrated for the $\beta$-cut.
		}
		\label{fig:size_comp_N_5_2}
	\end{figure} 
	
\newpage
	
\begin{figure}[H]
		\centering
		\begin{subfigure}{\textwidth}
			\centering
			\includegraphics[width=\textwidth]{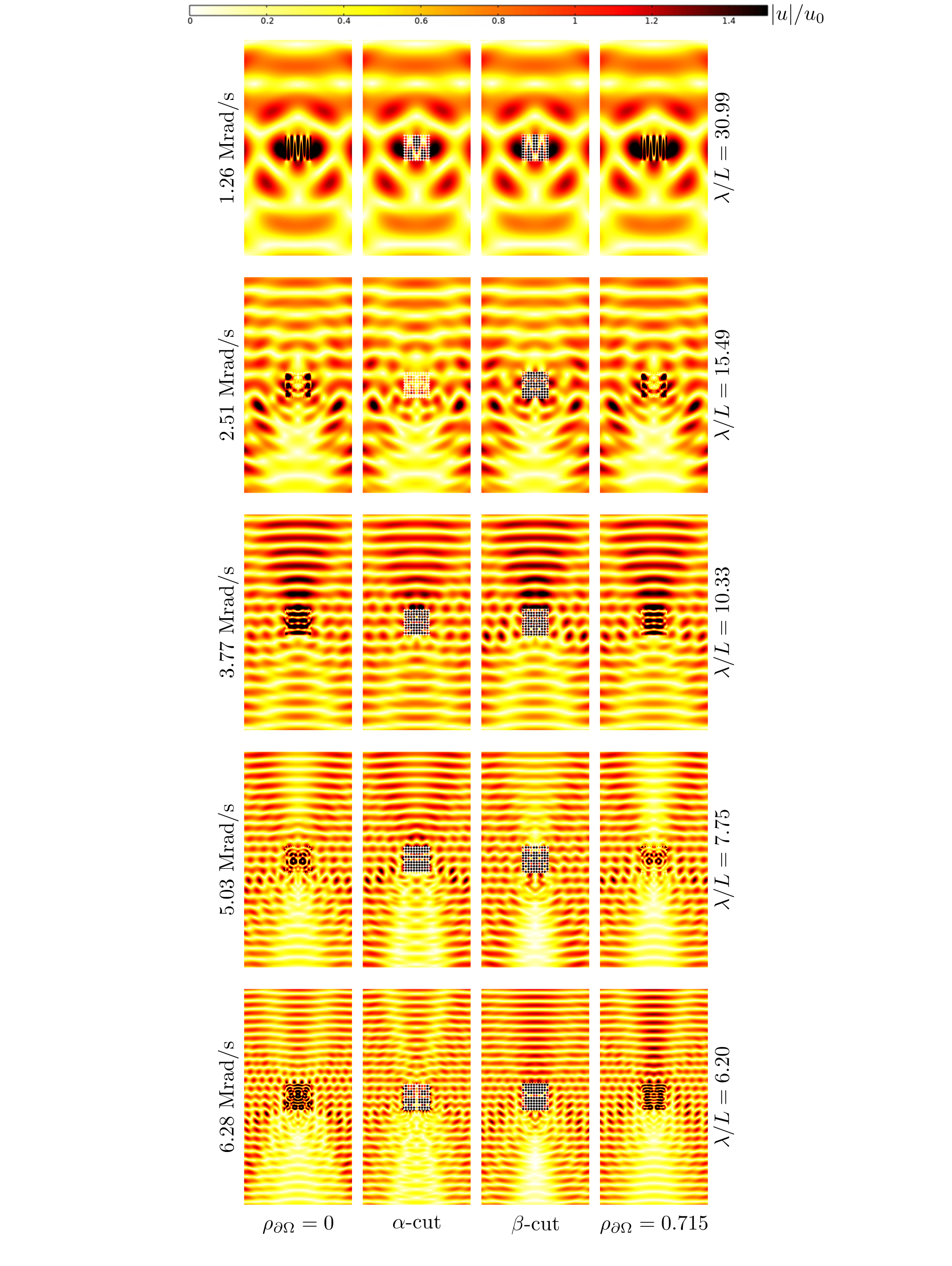}
		\end{subfigure}
		\caption{
		    Response at lower frequencies of a $10L\times10L$ metamaterial block as it interacts with a pressure plane wave that propagates vertically downwards.
		    The normalized displacement fields shown in each column correspond to the resulting scattering pattern for five types of implementation, from left to right:
		    \textit{($\rho_{\partial \Omega}=0$)} micromorphic simulations without interface contribution,
		    \textit{($\alpha$-cut)} microstructured simulation of the architected material using $\alpha$-type boundaries,
		    \textit{($\beta$-cut)} microstructured simulation of the architected material using $\beta$-type boundaries, and
		    \textit{($\rho_{\partial \Omega}=-0.715$)} enhanced micromorphic simulations using the surface inertia contribution calibrated for the $\beta$-cut.
		}
		\label{fig:size_comp_N_10_1}
	\end{figure}

\begin{figure}[H]
		\centering
		\begin{subfigure}{\textwidth}
			\centering
			\includegraphics[width=\textwidth]{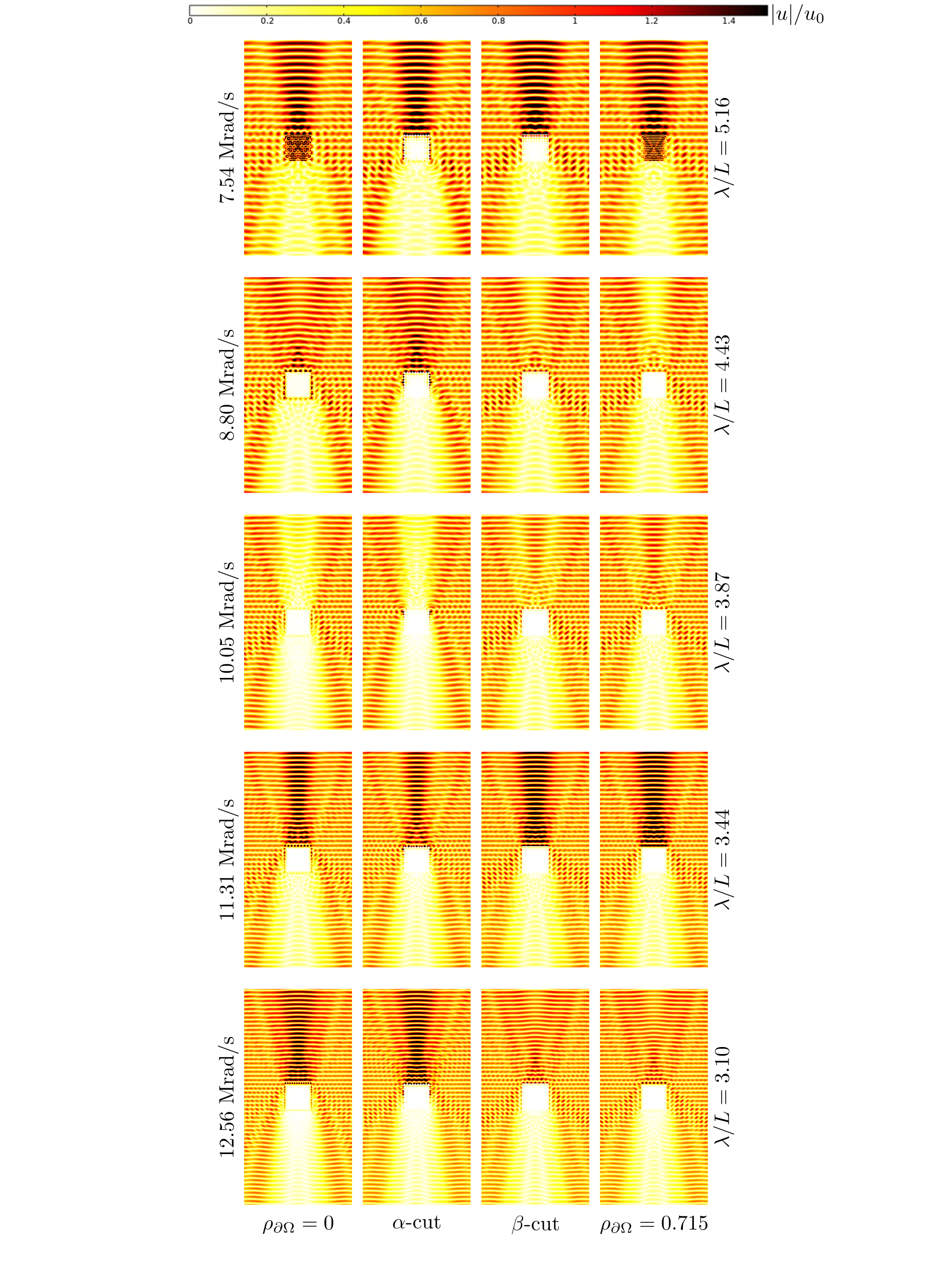}
		\end{subfigure}
		\caption{
		    Response at higher frequencies of a $10L\times10L$ metamaterial block as it interacts with a pressure plane wave that propagates vertically downwards.
		    The normalized displacement fields shown in each column correspond to the resulting scattering pattern for five types of implementation, from left to right:
		    \textit{($\rho_{\partial \Omega}=0$)} micromorphic simulations without interface contribution,
		    \textit{($\alpha$-cut)} microstructured simulation of the architected material using $\alpha$-type boundaries,
		    \textit{($\beta$-cut)} microstructured simulation of the architected material using $\beta$-type boundaries, and
		    \textit{($\rho_{\partial \Omega}=-0.715$)} enhanced micromorphic simulations using the surface inertia contribution calibrated for the $\beta$-cut.
		}
		\label{fig:size_comp_N_10_2}
	\end{figure}

\begin{figure}[H]
		\centering
		\begin{subfigure}{\textwidth}
			\centering
			\includegraphics[width=\textwidth]{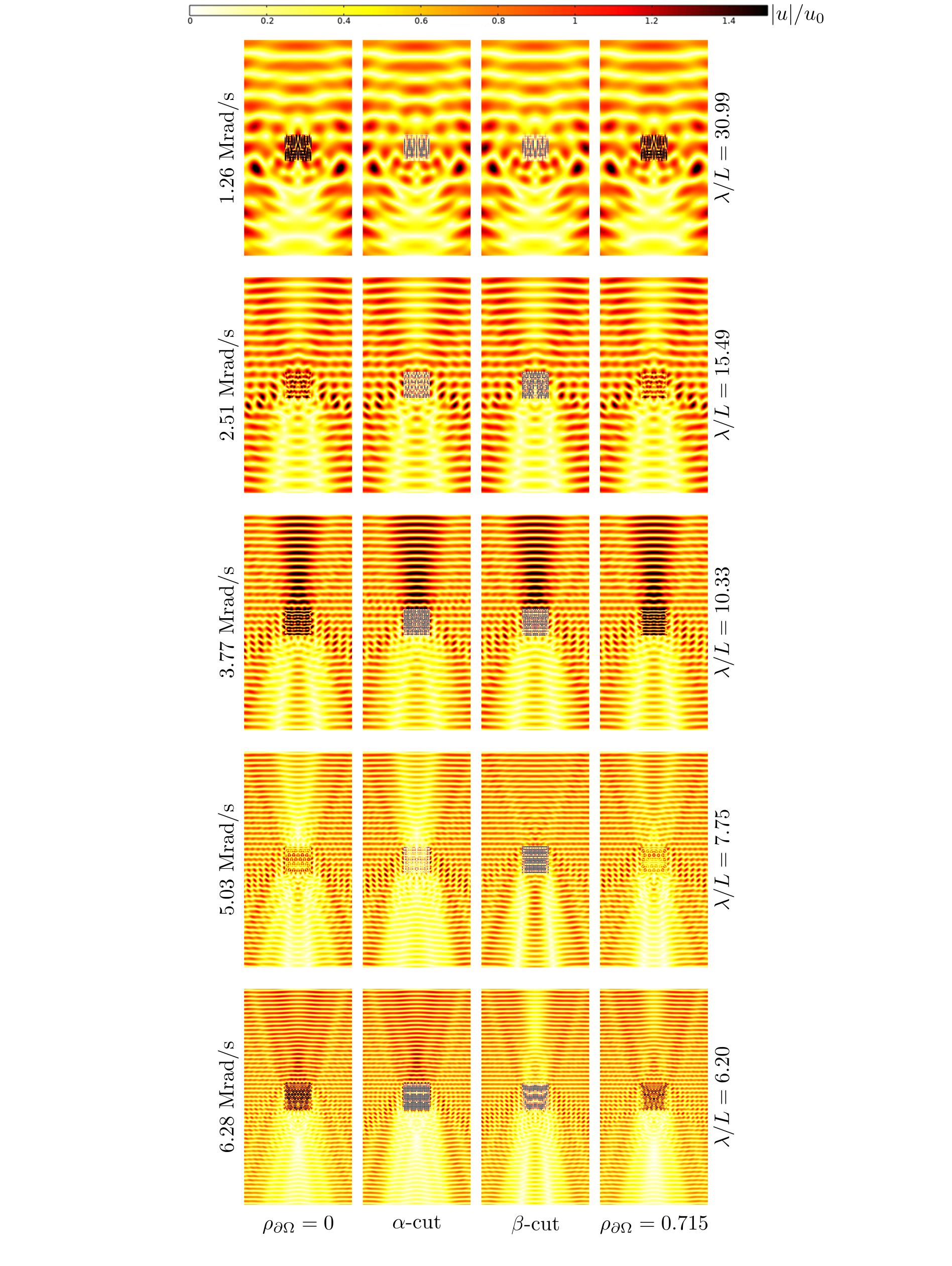}
		\end{subfigure}
		\caption{
		    Response at lower frequencies of a $20L\times20L$ metamaterial block as it interacts with a pressure plane wave that propagates vertically downwards.
		    The normalized displacement fields shown in each column correspond to the resulting scattering pattern for five types of implementation, from left to right:
		    \textit{($\rho_{\partial \Omega}=0$)} micromorphic simulations without interface contribution,
		    \textit{($\alpha$-cut)} microstructured simulation of the architected material using $\alpha$-type boundaries,
		    \textit{($\beta$-cut)} microstructured simulation of the architected material using $\beta$-type boundaries, and
		    \textit{($\rho_{\partial \Omega}=-0.715$)} enhanced micromorphic simulations using the surface inertia contribution calibrated for the $\beta$-cut.
		}
		\label{fig:size_comp_N_20_1}
	\end{figure}

\begin{figure}[H]
		\centering
		\begin{subfigure}{\textwidth}
			\centering
			\includegraphics[width=\textwidth]{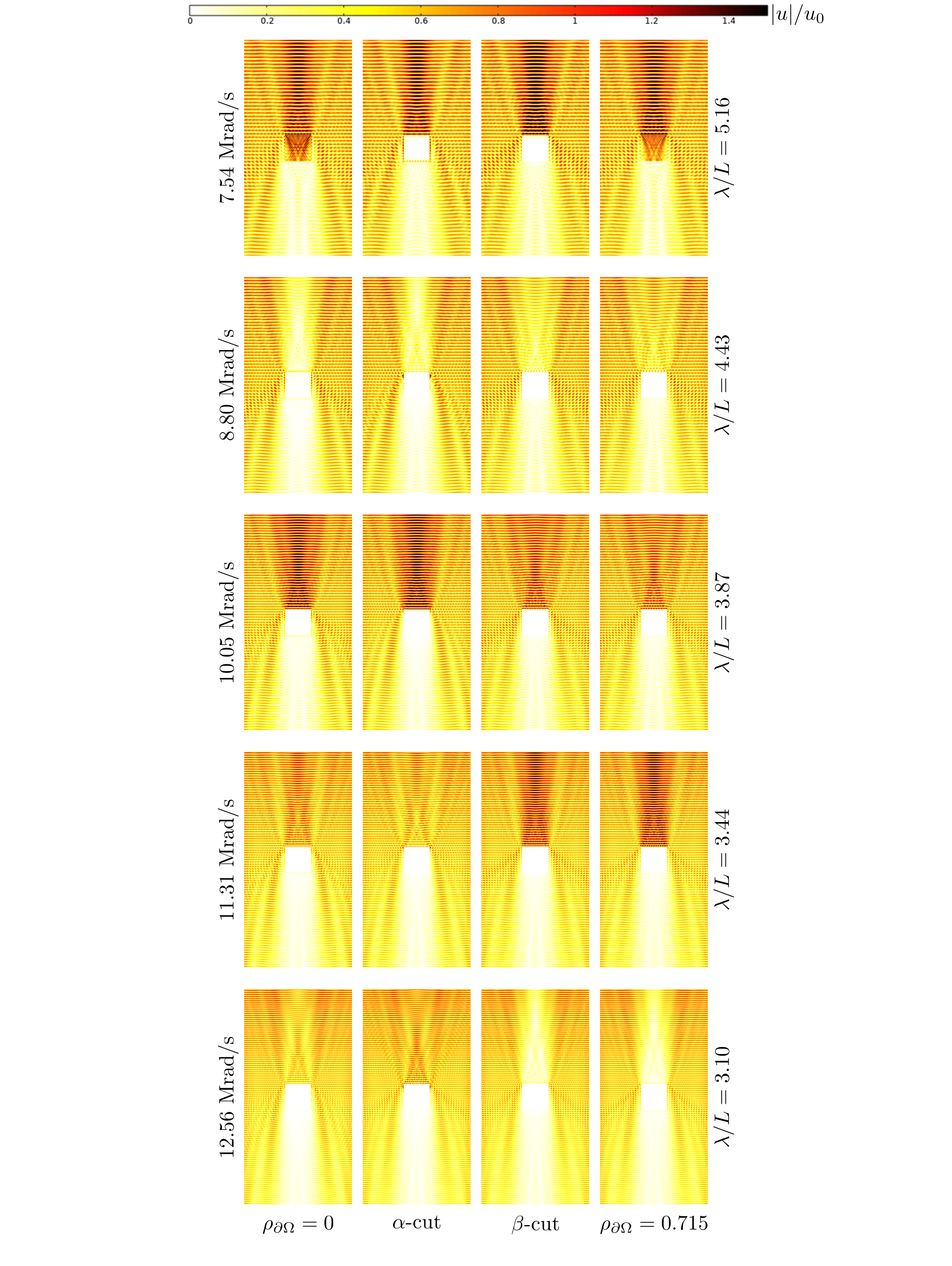}
		\end{subfigure}
		\caption{
		    Response at higher frequencies of a $20L\times20L$ metamaterial block as it interacts with a pressure plane wave that propagates vertically downwards.
		    The normalized displacement fields shown in each column correspond to the resulting scattering pattern for five types of implementation, from left to right:
		    \textit{($\rho_{\partial \Omega}=0$)} micromorphic simulations without interface contribution,
		    \textit{($\alpha$-cut)} microstructured simulation of the architected material using $\alpha$-type boundaries,
		    \textit{($\beta$-cut)} microstrucutred simulation of the architected material using $\beta$-type boundaries, and
		    \textit{($\rho_{\partial \Omega}=-0.715$)} enhanced micromorphic simulations using the surface inertia contribution calibrated for the $\beta$-cut.
		    }
		\label{fig:size_comp_N_20_2}
	\end{figure}

    As already remarked, we now proceed with the validation by analyzing the $\alpha$ and $\beta$ metamaterial's scattering patterns as well as their homogenized counterparts for $10L\times10L$ (Figs. \ref{fig:size_comp_N_10_1}-\ref{fig:size_comp_N_10_2}) and $20L\times20L$ (Figs. \ref{fig:size_comp_N_20_1}-\ref{fig:size_comp_N_20_2}) specimens' sizes.
    It can be directly inferred from Figs. \ref{fig:size_comp_N_10_1}-\ref{fig:size_comp_N_20_2} and tables \ref{tab:frequency-errors-N10}-\ref{tab:frequency-errors-N20} that the calculated value of $\rho_{\partial \Omega}=0.715$ continues allowing to retrieve the solution for the $\beta$-cut for all specimens' sizes and for a wide frequency range.
	
	\begin{table}[htbp]
      \centering
      \begin{threeparttable}
        \caption{Relative $L^2$ error [\%] of the normalized displacement field of the domain around the metamaterial block. The relative error of the micromorphic simulations is reported with respect to the $\alpha$- and $\beta$-cut types of boundaries of the corresponding microstructured simulation, at the prescribed angular frequencies for $N=10$.}
        \label{tab:frequency-errors-N10}
        \footnotesize
        \setlength{\tabcolsep}{2.7pt}
        \renewcommand{\arraystretch}{1.14}
        \begin{tabular}{@{}l S[table-format=-1.3]
          *{10}{S[table-format=2.2]} S[table-format=2.2]@{}}
          \freqheader
          $\alpha$ &  0     & 5.47 & 21.29 & 16.44 & 27.32 & 24.82 & 21.95 & 20.77 & 19.89 & 18.72 & 14.72 & 19.14 \\
          $\alpha$ & 0.715 & 7.54 & 22.39 & 13.35 & 33.91 & 32.37 & 31.69 & 33.81 & 35.77 & 32.34 & 34.67 & 27.79 \\
          \addlinespace[2pt]
          $\beta$  &  0     & 6.37 & 16.07 & 13.87 & 31.44 & 36.30 & 21.67 & 27.08 & 30.98 & 25.12 & 31.86 & 24.08 \\
          $\beta$  & 0.715 & 9.88 & 19.34 & 11.30 & 24.27 & 22.45 & 10.05 &  9.95 &  9.16 &  6.72 &  9.12 & 13.22 \\
          \bottomrule
        \end{tabular}
      \end{threeparttable}
    \end{table}

\newpage	
	
	\begin{table}[htbp]
      \centering
      \begin{threeparttable}
        \caption{Relative $L^2$ error [\%] of the normalized displacement field of the domain around the metamaterial block. The relative error of the micromorphic simulations is reported with respect to the $\alpha$- and $\beta$-cut types of boundaries of the corresponding microstructured simulation, at the prescribed angular frequencies for $N=20$.}
        \label{tab:frequency-errors-N20}
        \footnotesize
        \setlength{\tabcolsep}{2.7pt}
        \renewcommand{\arraystretch}{1.14}
        \begin{tabular}{@{}l S[table-format=-1.3]
          *{10}{S[table-format=2.2]} S[table-format=2.2]@{}}
          \freqheader
          $\alpha$ &  0     & 13.17 & 10.77 & 18.73 & 19.91 & 16.37 & 15.55 & 15.89 & 12.43 & 14.76 & 12.28 & 14.99 \\
          $\alpha$ & 0.715 & 11.58 & 10.31 & 13.64 & 30.06 & 25.34 & 26.52 & 29.10 & 27.88 & 33.41 & 34.41 & 24.23 \\
          \addlinespace[2pt]
          $\beta$  &  0     & 11.55 & 15.26 & 15.07 & 34.05 & 29.30 & 18.18 & 22.71 & 21.73 & 23.08 & 29.44 & 22.04 \\
          $\beta$  & 0.715 & 12.71 & 11.10 &  8.04 & 29.78 & 19.70 &  8.65 & 10.22 &  6.96 &  6.26 &  7.36 & 12.08 \\
          \bottomrule
        \end{tabular}
      \end{threeparttable}
    \end{table}

	\subsection{Correlation Between Truncation Plane and Surface Density: Validation of the Interface Inertia Concept on the $\gamma$-Cut}
	
	We now consider the $\gamma$-cut as an additional and independent validation case for the proposed interface-inertia concept.
	This cut is generated by a different truncation plane of the same periodic cross-cell microstructure and therefore provides a useful mean to confirm the new interface inertia concept.
	
	Using the same calibration strategy adopted for the other interfaces, the effective surface density associated with the $\gamma$-cut is found to be $\evalat{\surfdens}{\gamma}=0.3575$ $\rm kg/m^2$.
	
\begin{figure}[H]
		\centering
		\begin{subfigure}{\textwidth}
			\centering
			\includegraphics[width=0.9\textwidth]{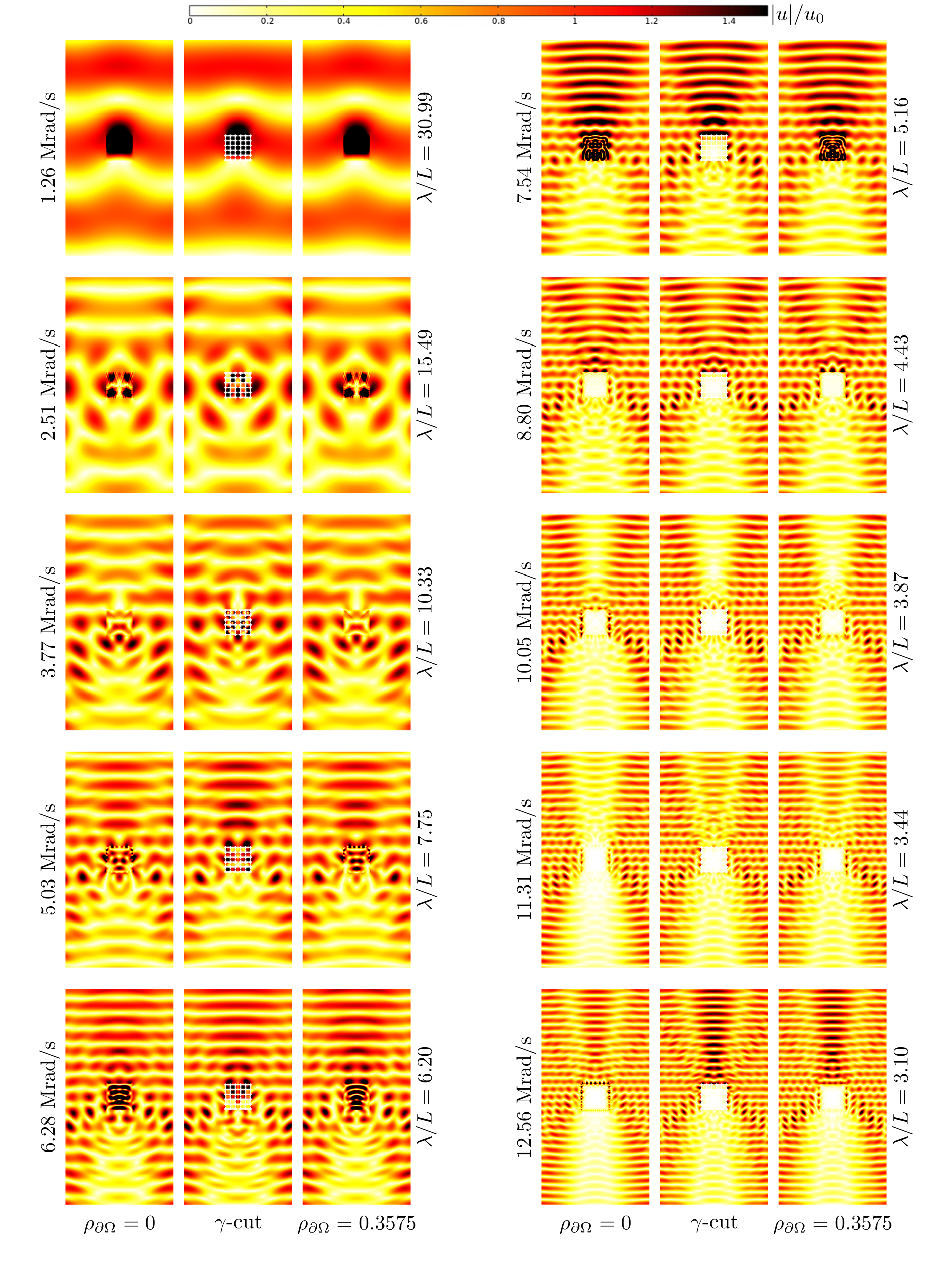}
		\end{subfigure}
		\caption{
			Study of the response of a $5.75L\times5.75L$ metamaterial block as it interacts with a pressure plane wave that propagates vertically downwards.
			The normalized displacement fields shown in each column correspond to the resulting scattering pattern for five types of implementation:
			\textit{($\rho_{\partial \Omega}=0$)} micromorphic simulations with coherent interfaces,
			\textit{($\delta-$cut)} direct simulation of the architected material using $\delta$-type boundaries, and
			\textit{($\rho_{\partial \Omega}=-0.3575$)} enhanced micromorphic simulations using the surface inertia contribution targeted for the $\delta-$cut.
		}
		\label{fig:delta_comp_N_5.75}
	\end{figure}

\begin{figure}[H]
		\centering
		\begin{subfigure}{\textwidth}
			\centering
			\includegraphics[width=0.9\textwidth]{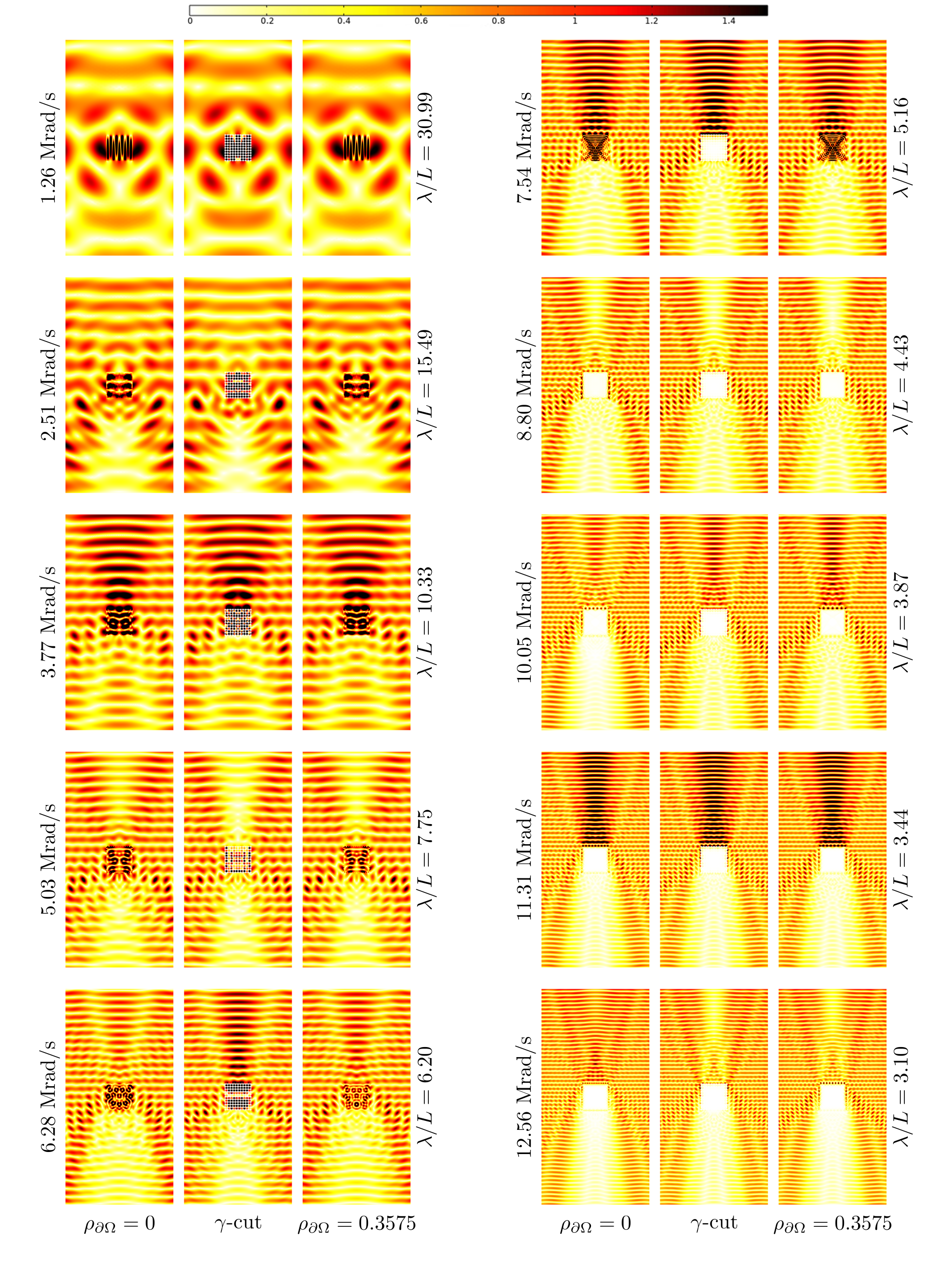}
		\end{subfigure}
		\caption{
			Study of the response of a $10.75L\times10.75L$ metamaterial block as it interacts with a pressure plane wave that propagates vertically downwards.
			The normalized displacement fields shown in each column correspond to the resulting scattering pattern for five types of implementation:
			\textit{($\rho_{\partial \Omega}=0$)} micromorphic simulations with coherent interfaces,
			\textit{($\delta-$cut)} direct simulation of the architected material using $\delta$-type boundaries, and
			\textit{($\rho_{\partial \Omega}=-0.3575$)} enhanced micromorphic simulations using the surface inertia contribution targeted for the $\delta-$cut.
		}
		\label{fig:delta_comp_N_10.75}
	\end{figure}

\begin{figure}[H]
		\centering
		\begin{subfigure}{\textwidth}
			\centering
			\includegraphics[width=0.9\textwidth]{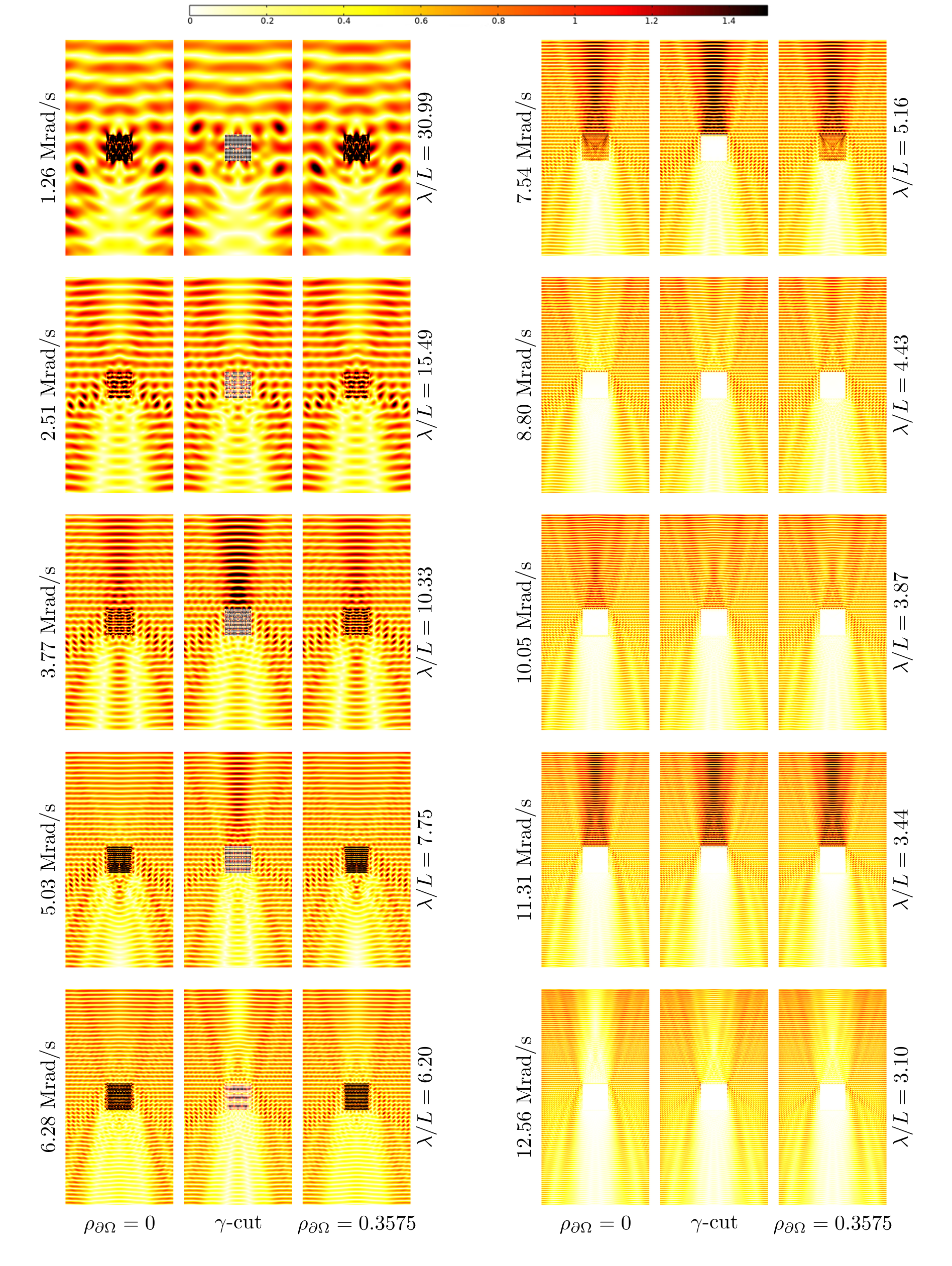}
		\end{subfigure}
		\caption{
			Study of the response of a $20.75L\times20.75L$ metamaterial block as it interacts with a pressure plane wave that propagates vertically downwards.
			The normalized displacement fields shown in each column correspond to the resulting scattering pattern for five types of implementation:
			\textit{($\rho_{\partial \Omega}=0$)} micromorphic simulations with coherent interfaces,
			\textit{($\delta-$cut)} direct simulation of the architected material using $\delta$-type boundaries, and
			\textit{($\rho_{\partial \Omega}=-0.3575$)} enhanced micromorphic simulations using the surface inertia contribution targeted for the $\delta-$cut.
		}
		\label{fig:delta_comp_N_20.75}
	\end{figure}

	The validation results in Figs \ref{fig:delta_comp_N_5.75}-\ref{fig:delta_comp_N_20.75} show that this value allows the enhanced relaxed micromorphic model to reproduce the main features of the fully resolved $\gamma$-cut response for a large frequency range and for different specimen’s sizes.
	
	The agreement of the relaxed micromorphic simulations with the fully microstructured ones improves both in terms of the overall structure of the scattered field and in terms of the displacement distribution around the metamaterial block.
	
	The same value of $\evalat{\surfdens}{\gamma}$ is used for all the $\gamma$-cut specimens considered.
	The improvement	obtained for different specimen sizes shows that the calibrated surface density is not tied to a particular finite sample, but is instead associated with the geometry of the $\gamma$-type interface.
	This supports the interpretation of rho as an effective homogenized parameter encoding the additional inertia generated by the boundary truncation.
	
	The $\gamma$-cut validation therefore confirms that different truncation planes of the same metamaterial may require different effective interface inertias, even when the bulk relaxed micromorphic parameters are kept unchanged.
	The proposed surface-inertia term enables the homogenized model to account for this truncation-dependent boundary effect.

    \begin{figure}[htbp]
		\centering
		\begin{subfigure}{\textwidth}
			\centering
			\includegraphics[width=\textwidth]{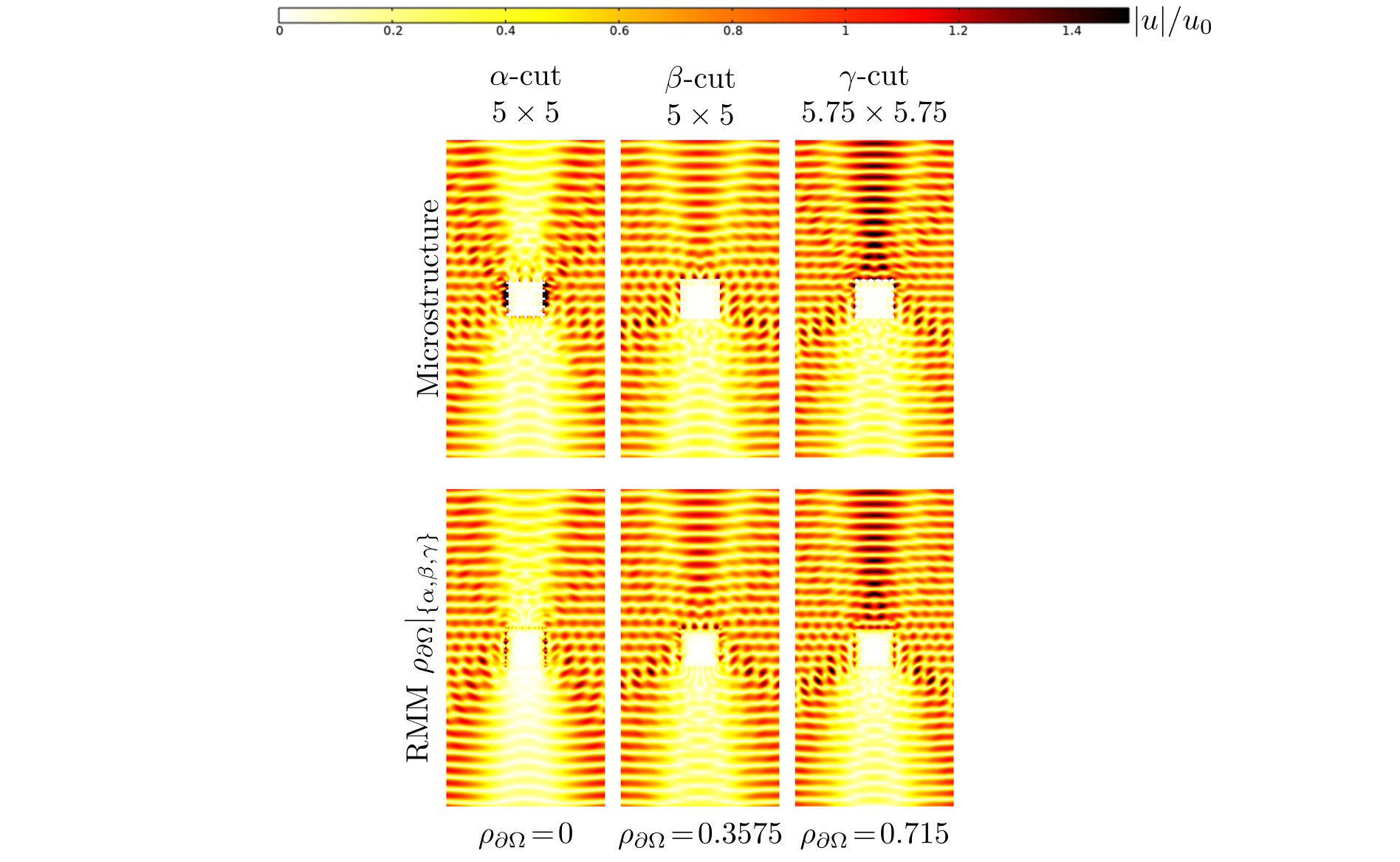}
		\end{subfigure}
		\caption{
		    Normalized displacement field of a pressure plane wave as it propagates at $12.56$ M\,rad/s vertically downwards and interacts with a metamaterial block. 
		    Each column shows the scattering patterns for four types of metamaterial-boundaries and the size of each block is given as a multiple of $L$.
		    The rows correspond to different simulation approaches:
		    \textit{(Microstructure)} shows the results of the direct simulation of the architected material, and
		    \textit{(RMM $\rho_{\partial \Omega}$)} present the enhanced micromorphic solution by using surface inertia $\evalat{\surfdens}{\beta}=0$, $\evalat{\surfdens}{\beta}=0.715$ and $\evalat{\surfdens}{\gamma}=0.
		    3575$.
		}
		\label{fig:alpha_beta_gamma_delta-comparison}
	\end{figure}

	\subsubsection{Frequency Analysis of the $\gamma$-Cut}
	
	\begin{table}[H]
		\centering
		\begin{threeparttable}
			\caption{Relative $L^2$ error [\%] of the normalized displacement field of the domain around the metamaterial block. The relative error of the micromorphic simulations is reported respect to microstructured simulations of a block with $\gamma-$cut types of boundaries and the corresponding reference size, at the prescribed angular frequencies.}
			\label{tab:frequency-errors-delta}
			\footnotesize
			\setlength{\tabcolsep}{2.15pt}
			\renewcommand{\arraystretch}{1.14}
			\begin{tabular}{@{}S[table-format=2.3] S[table-format=-1.5]
					*{10}{S[table-format=2.2]} S[table-format=2.2]@{}}
				\freqheader
				5.75 &  0       &  4.75 & 14.19 & 17.56 & 19.25 & 20.16 & 16.68 & 20.74 & 26.53 & 27.53 & 27.69 & 19.51 \\
				5.75 & 0.3575 &  5.39 & 10.86 & 13.80 & 11.66 & 14.66 & 11.69 &  9.59 & 11.15 & 13.27 & 15.80 & 11.79 \\
				\addlinespace[2pt]
				10.75 &  0       & 16.64 & 24.42 & 11.24 & 26.38 & 28.75 & 17.01 & 20.71 & 22.14 & 16.77 & 24.42 & 20.85 \\
				10.75 & 0.3575 & 18.27 & 25.11 &  6.90 & 15.73 & 18.51 & 10.21 &  8.91 &  9.22 &  9.33 & 15.32 & 13.75 \\
				\addlinespace[2pt]
				20.75 &  0       & 21.13 & 11.72 & 21.36 & 29.29 & 24.60 & 13.46 & 15.95 & 15.69 & 14.06 & 19.66 & 18.69 \\
				20.75 & 0.3575 & 20.17 &  9.71 & 18.78 & 26.68 & 19.88 & 10.70 &  8.66 &  5.64 &  4.69 & 10.85 & 13.58 \\
				\bottomrule
			\end{tabular}
		\end{threeparttable}
	\end{table}
	
\newpage 
	
\FloatBarrier	

\section{Conclusion and Perspectives}
	
	In this work, we have introduced an interface-inertia-enhanced relaxed micromorphic model for the homogenized description of finite-size mechanical metamaterials.
	The central motivation was the observation that, when the wavelength becomes comparable to the unit-cell size, the classical separation of scales progressively breaks down and the dynamic response of a finite metamaterial specimen becomes sensitive to the way in which the underlying microstructure is truncated at the external boundary.
	
	Fully resolved simulations showed that different cuts of the same periodic lattice can lead to distinct scattering responses, even though the bulk microstructure is identical.
	This demonstrates that the boundary of a finite metamaterial cannot always be treated as a purely geometrical surface carrying only standard tractions.
	Instead, at sufficiently high frequencies, the interface behaves as an effective dynamic region whose inertial contribution depends on a boundary-layer mass distribution generated by the truncation plane.
	
	To account for this effect within a homogenized framework, we enriched the relaxed micromorphic model by adding a kinetic surface energy at the boundary of the metamaterial block.
	This term leads, through a variationally consistent derivation, to an additional inertial contribution in the boundary condition.
	The resulting model preserves the bulk relaxed micromorphic structure while allowing the interface to carry an effective surface inertia.
	In this way, specimens with the same bulk parameters but different boundary truncations can be distinguished at the homogenized level.
	
	The numerical results confirm the relevance of this novel approach.
	The standard relaxed micromorphic model with coherent interface conditions accurately describes the reference $\alpha$-cut response, while 	different truncations require non-zero effective interface densities. 
	Once the calibrated surface density is introduced, the enhanced model reproduces the main features of the corresponding fully resolved scattering fields over a broad frequency range and for different specimen sizes.
	This shows that the identified interface density is not merely a sample-dependent fitting parameter, but a homogenized measure of the truncation-dependent boundary inertia. 
	
	The validation on several cuts further supports the physical interpretation of the proposed surface inertia.
	Different truncation planes give rise to different calibrated values of the interface density, consistently with the fact that they generate different interfacial boundary-layer mass distributions.
	The model is therefore able to encode boundary effects that are completely absent from classical bulk-only homogenized formulations.
	
	Some small discrepancies remain between the enhanced relaxed micromorphic predictions and the fully resolved simulations, especially in intermediate frequency regimes where the wavelength is close to the macroscopic specimen size and the interaction between bulk dispersion, boundary scattering, and local microstructural resonances becomes more complex.
	These residual differences are not unexpected, since the present model introduces only a scalar surface-inertia correction.
	The	remaining mismatch suggests that additional interface mechanisms may contribute to the boundary response.
	
	A first natural extension of the present formulation is the introduction of a non-isotropic interface inertia.
	In the present work, the surface density is described by a scalar parameter, which assumes that the inertial correction is the same in all displacement directions.
	However, the boundary microstructure generated by a given truncation plane may have a directional mass distribution.
	In that case, the interface inertia should be represented by a surface inertia tensor rather than by a scalar density.
	Such an extension would allow the model to distinguish between normal and tangential inertial effects at the boundary and could further improve the agreement with fully resolved simulations, especially when mode conversion or mixed pressure-shear scattering is significant.
	
	A second possible extension is the introduction of a surface elastic energy in addition to the kinetic surface energy proposed here.
	Such a term would allow the interface not only to carry inertia, but also to store elastic energy associated with localized deformation mechanisms at the boundary.
	This could improve the description of interface-induced scattering, especially in regimes where the boundary response is governed not only by mass redistribution but also by local stiffness effects.
	
	In future works, we will investigate these extensions by exploring orthotropic surface inertia and surface elasticity into the relaxed micromorphic variational formulation.
	
	Overall, the present results show that interface inertia provides a powerful and physically motivated correction to homogenized metamaterial models when scale separation is no longer valid.
	The proposed formulation opens the way toward predictive continuum models for finite metamaterial building blocks, where bulk dispersion and truncation-dependent boundary dynamics can be treated within a single variational framework.
	This is an essential step toward the efficient simulation of large structures assembled from finite metamaterial components, for which fully resolved microstructured computations would be prohibitively expensive.

\section*{Acknowledgement} 
The authors acknowledge support from the European Commission through the funding of the ERC Consolidator Grant META-LEGO, no. 101001759. 

\clearpage
	
	\appendix

\section{Differentiated Response of a $20L\times20L$ Metamaterial Block}
	
	In this appendix, results analogous to those presented in Section 2.1 are reported for a bigger specimen size (20 × 20 unit cells). We remark that truncation-related differences continue to be important for sufficiently small wavelengths.
	
	\begin{figure}[H]
		\begin{subfigure}{\textwidth}
			\centering
			\includegraphics[width=\textwidth]{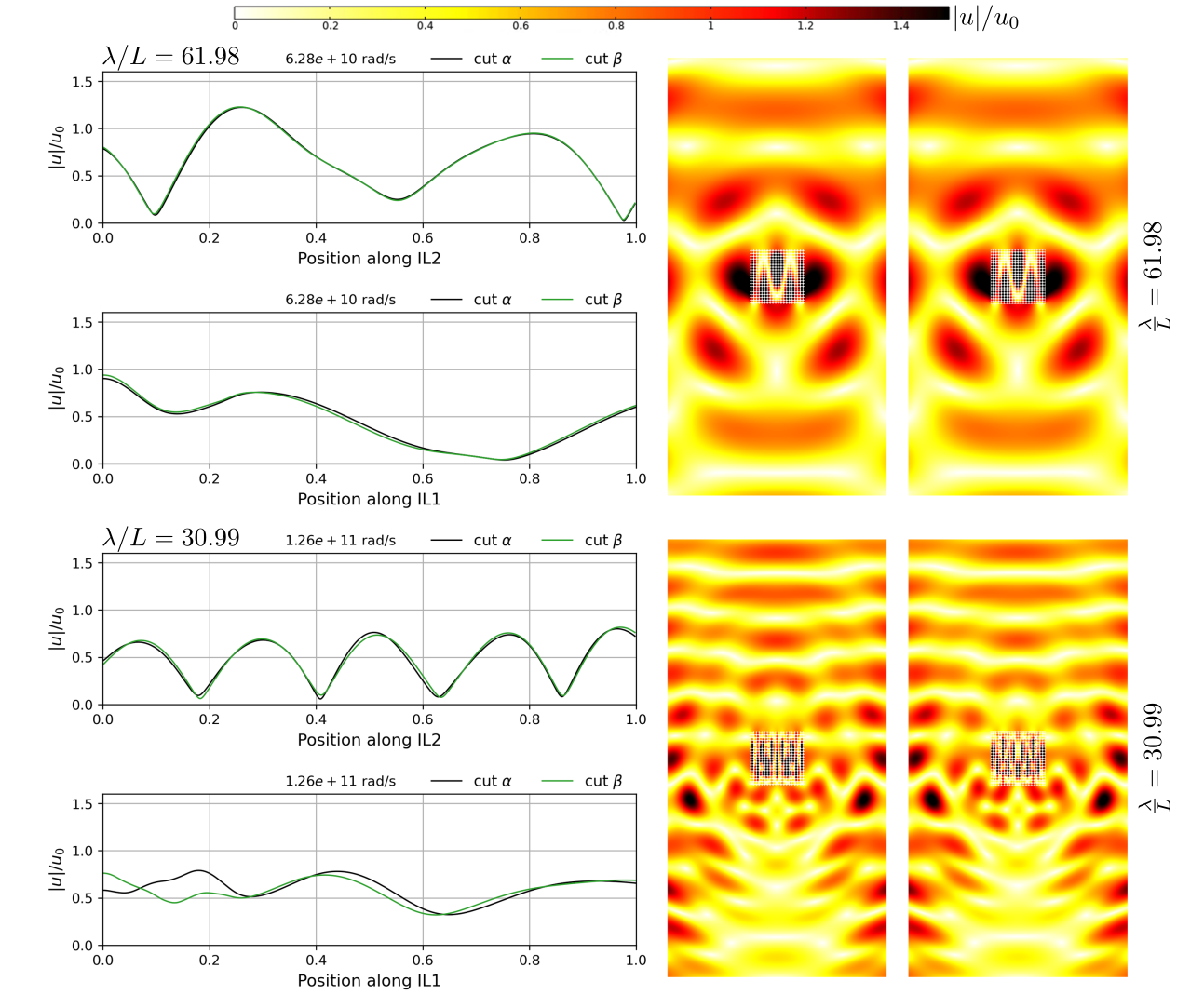}
		\end{subfigure}
		\caption{
		    Differentiated response of the scattering pattern, driven by the boundary contribution of $\alpha$ and $\beta$ interfaces, as a pressure wave interacts with a $20L\times20L$ metamaterial block.
		    The normalized displacement along inspection lines 1 and 2, as well as the normalized displacement field associated with a metamaterial block with $\alpha$- and $\beta$-type boundaries, are shown for the frequencies:
			\textit{(upper)} $0.63\ \mathrm{M\,rad/s}$ and
			\textit{(lower)} $1.26\ \mathrm{M\,rad/s}$.
			The corresponding wavelengths are much bigger than the size of the unit cell: no significant differences between $\alpha$- and $\beta$-cuts can be appreciated.
			The hypothesis of scale separation can be expected to hold true.
		}
		\label{fig:mstd_differences-20_low}
	\end{figure}

	\begin{figure}[H]
		\begin{subfigure}{\textwidth}
			\centering
			\includegraphics[width=\textwidth]{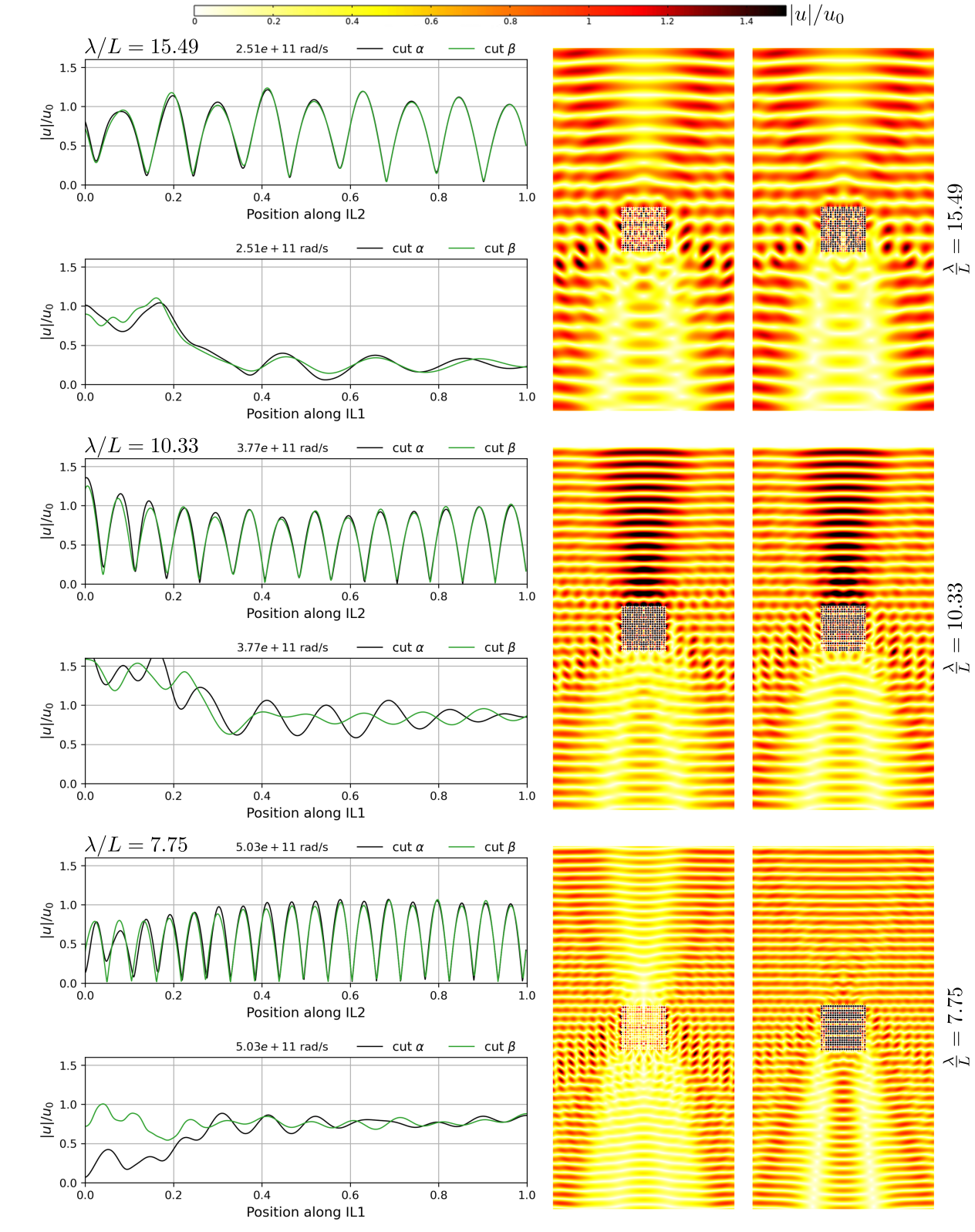}
		\end{subfigure}
		\caption{
		    Differentiated response of the scattering pattern, driven by the boundary contribution of $\alpha$ and $\beta$ interfaces, as a pressure wave interacts with a $20L\times20L$ metamaterial block.
		    The normalized displacement along inspection lines 1 and 2, as well as the normalized displacement field associated with a metamaterial block with $\alpha$- and $\beta$-type boundaries, are shown for the frequencies:
			\textit{(upper)} $2.51\ \mathrm{M\,rad/s}$,
			\textit{(middle)} $3.77\ \mathrm{M\,rad/s}$,
			and \textit{(lower)} $5.03\ \mathrm{M\,rad/s}$.
			The corresponding wavelengths become comparable to the size of the unit cell: small differences between $\alpha$- and $\beta$-cut start being observable.
			The wavelength $\lambda = 10.33L$ is the one at which the first differences between the two cuts start becoming apparent: we will then suppose that at this threshold value and below, separation of scales is not expected to hold true anymore.
		}
		\label{fig:mstd_differences-20_mid}
	\end{figure}

	\begin{figure}[H]
		\begin{subfigure}{\textwidth}
			\centering
			\includegraphics[width=\textwidth]{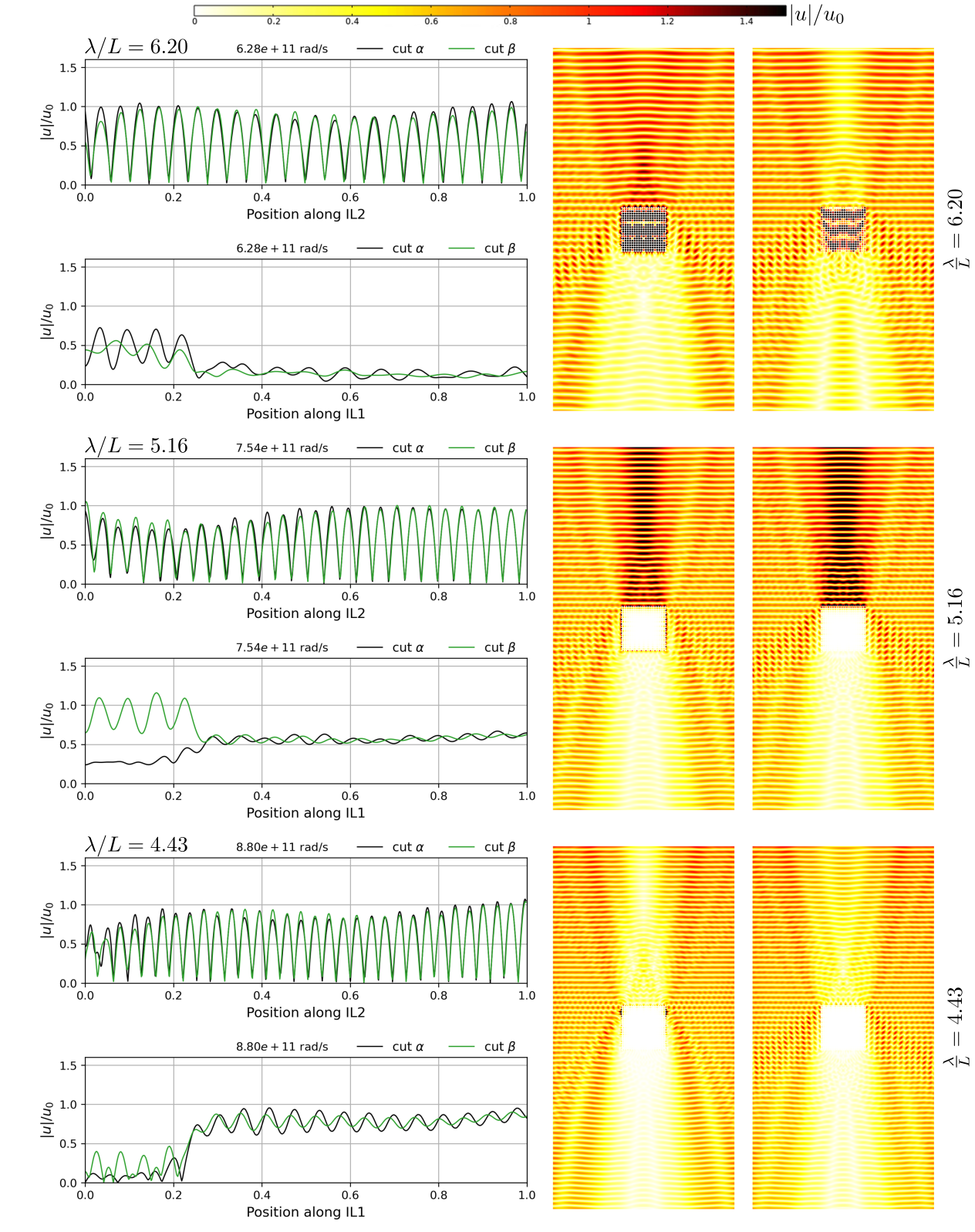}
		\end{subfigure}
		\caption{
		    Differentiated response of the scattering pattern, driven by the boundary contribution of $\alpha$ and $\beta$ interfaces, as a pressure wave interacts with a $20L\times20L$ metamaterial block.
		    The normalized displacement along inspection lines 1 and 2, as well as the normalized displacement field associated with a metamaterial block with $\alpha$- and $\beta$-type boundaries, are shown for the frequencies:
			\textit{(upper)} $6.28\ \mathrm{M\,rad/s}$,
			\textit{(middle)} $7.54\ \mathrm{M\,rad/s}$,
			and \textit{(lower)} $8.80\ \mathrm{M\,rad/s}$.
			The corresponding wavelengths become closer to the unit cell's size: larger differences between $\alpha$- and $\beta$-cuts can be observed. Separation of scales hypothesis breaks down.
		}
		\label{fig:mstd_differences-20_high}
	\end{figure}

	\begin{figure}[H]
		\begin{subfigure}{\textwidth}
			\centering
			\includegraphics[width=\textwidth]{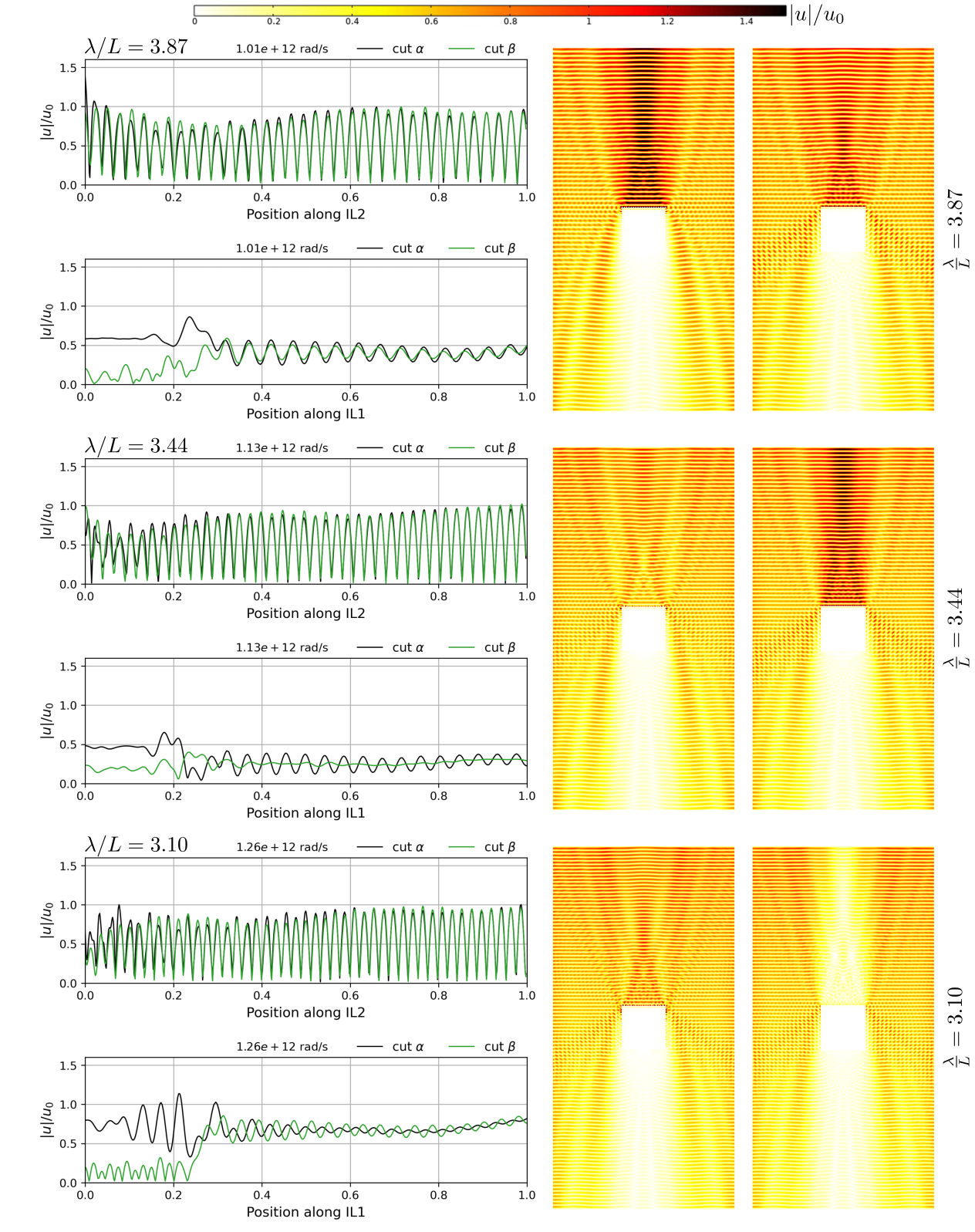}
		\end{subfigure}
		\caption{
		    Differentiated response of the scattering pattern, driven by the boundary contribution of $\alpha$ and $\beta$ interfaces, as a pressure wave interacts with a $20L\times20L$ metamaterial block.
		    The normalized displacement along inspection lines 1 and 2, as well as the normalized displacement field associated with a metamaterial block with $\alpha$- and $\beta$-type boundaries, are shown for the frequencies:
			\textit{(upper)} $10.05\ \mathrm{M\,rad/s}$,
			\textit{(middle)} $11.31\ \mathrm{M\,rad/s}$,
			and \textit{(lower)} $12.56\ \mathrm{M\,rad/s}$.
			The corresponding wavelengths are very close to the unit cell's size: the differences between $\alpha$- and $\beta$-cuts are important. Separation of scales hypothesis does not hold true anymore.
		}
		\label{fig:mstd_differences-20_band_gap}
	\end{figure}

{\footnotesize
\printbibliography
}

\end{document}